\documentclass[11pt]{article}
\font\eufm=eufm10
\def\frak#1{\hbox{\eufm#1}}
\newcommand{\bd}{
\begin{document}}
\newcommand{\ed}{\end{document}}
\newcommand{\be}{\begin{enumerate}}
\newcommand{\ee}{\end{enumerate}}
\newcommand{\bi}{\begin{itemize}}
\newcommand{\ei}{\end{itemize}}
\newcommand{\ba}{\begin{array}}
\newcommand{\ea}{\end{array}}
\newcommand{\vs}{\vspace*{0.3\baselineskip}}
\newcommand{\vsm}{\vspace*{-0.3\baselineskip}}
\newtheorem{defi}{Definition}[section]
\newtheorem{tw}[defi]{Theorem}
\newtheorem{prop}[defi]{Proposition}
\newtheorem{lem}[defi]{Lemma}
\newtheorem{re}[defi]{Remark}
\newtheorem{col}[defi]{Corollary}
\newtheorem{ex}[defi]{Examples}
%
\newcommand{\Om}{\Omega}
\newcommand{\om}{\omega}
\newcommand{\G}{\Gamma}
\newcommand{\D}{\Delta}
\renewcommand{\d}{\delta}
\newcommand{\ga}{\gamma}
\newcommand{\eps}{\epsilon}
\newcommand{\C}{C}
\newcommand{\R}{R}
%
\newcommand{\ove}{\overline}
\newcommand{\ms}{\oplus}
\newcommand{\mt}{\otimes}
\newcommand{\dz}{\wedge}
\newcommand{\lra}{\longrightarrow}
\newcommand{\ti}{\tilde}
\newcommand{\rel}{\mbox{$\,$\rule[0.5ex]{1.1em}{0.2pt}$\triangleright\,$}}
\newcommand{\dow}{\hspace*{\fill}\rule{1.6ex}{1.6ex}\hspace*{1em}}  
\newcommand{\dowl}{\hspace*{\fill}\rule{1ex}{1ex}\hspace*{1em}} 
\newcommand{\sd}{\hspace{0.3ex}\tiny{\rhd\mbox{\hspace{-2ex}}<}\hspace{0.3ex}}
%
\newcommand{\g}{\frak g}
\newcommand{\ab}{\frak a}
\newcommand{\bb}{\frak b}
\newcommand{\h}{\frak h}
\newcommand{\got}{\frak t}
\newcommand{\hd}{\hat{\d}}
\newcommand{\oml}{\Omega_L^{1/2}}
\newcommand{\omr}{\Omega_R^{1/2}}
\newcommand{\omh}{\Omega^{1/2}}
\newcommand{\lo}{\lambda_0}
\newcommand{\ro}{\rho_0}
\newcommand{\sA}{\mbox{$\cal A\,$}}
\newcommand{\lma}{\Lambda^{max}}
\newcommand{\timh}{\times_h}
\newcommand{\Gd}{\G^{(2)}}
\newcommand{\el}{e_L}
\newcommand{\er}{e_R}
\newcommand{\GG}{\G_1\times\G_2}
\newcommand{\gdot}{\hspace{-0.1em}\cdot\hspace{-0.1em}}
\newcommand{\tran}{\frown\hspace{-2.2ex}|\hspace{1.9ex}}
%
\newcommand{\al}{a_L}
\newcommand{\ar}{a_R}
\newcommand{\bl}{b_L}
\newcommand{\br}{b_R}
\newcommand{\adab}{\varphi} 
\newcommand{\adbb}{\psi} 
\newcommand{\adaa}[1]{\mbox{$|\det Ad #1 |_{\ab}|^{1/2}$ }}
\newcommand{\pra}{\mbox{$P_A$}}
\newcommand{\prb}{\mbox{$P_B$}}
\newcommand{\lmA}{\mbox{$\lma T_eA$ }}
\newcommand{\lmB}{\mbox{$\lma T_eB$ }}
\newcommand{\Adaa}{\mbox{ $\det Ad(a)|_{\ab}$ }}
\renewcommand{\Adaa}{\psi_A}
\newcommand{\Adab}{\mbox{ $\det Ad(a)|_{\bb}$ }}
\renewcommand{\Adab}{\psi_B}
\newcommand{\Adba}{\mbox{ $\det Ad(b)|_{\ab}$ }}
\renewcommand{\Adba}{\varphi_A}
\newcommand{\Adbb}{\mbox{ $\det Ad(b)|_{\bb}$ }}
\renewcommand{\Adbb}{\varphi_B}
%
\newcommand{\la}[2]{\Lambda_{#1#2}}
\newcommand{\kad}{ad^{\#}}
\newcommand{\wl}[1]{\vphantom{X}_{#1}{\G}}
\newcommand{\te}{\tilde{e}}
\newcommand{\sD}{\mbox{$\cal D\,$}}
\newcommand{\notka}[1]{}

\title{From Double Lie Groups to Quantum Groups.}

\author{Piotr Stachura\thanks{Department of Mathematical Methods in Physics,
University of Warsaw, Ho\.{z}a 74, 00-682 Warszawa, Poland}
\vspace{1ex}\thanks{The work was supported by Polish KBN grant No. 2 P03A 030 14 and
The Foundation for Polish Science.}}

\bd
\maketitle
\noindent{\bf Abstract\hspace{1em}} {\small It is shown that there is 
a $C^*$-algebraic quantum group related to  any double Lie group. 
An  algebra underlying this  quantum group is an algebra of a differential 
groupoid naturally associated with a double Lie group.}

\setcounter{section}{-1}
\section{Introduction}

The theory of objects we now call Quantum Groups is  from its very 
beginning related to group factorization. First constructions are
due to  Kac \cite{Kac}. Then an important step (in the context of
Kac-von Neumann algebras) was made by Majid \cite{maj1,maj2}. Afterwards,
Baaj and Skandalis created their theory of multiplicative unitaries 
\cite{BajSka}. From this moment on one  could expect that there is a $C^*$-algebraic
quantum group behind any (locally compact) group factorization. And indeed,
recently a very general construction of $C^*$-algebraic quantum groups
related to a group factorization was presented in \cite{Vaes-Vain}.

This paper is devoted to a construction of quantum groups from factorizations 
known as Double Lie Groups (DLG) \cite{WeLu}. In contrast to \cite{Vaes-Vain} 
we are able to use different, more geometric methods, because we deal with Lie groups.
This fact enables us to use  their rich geometric structure.
 
Since we are dealing with Lie groups, the objects used in our constructions
are defined at a smooth level. All of them have a clear geometric
interpretation as objects related to  differential groupoids associated 
naturally with DLGs. Starting  from a concrete algebra of 
such objects we  construct at this smooth level various mappings, 
which appear in the theory of quantum groups. Then by using density arguments 
we lift those mappings to a $C^*$-algebraic level.

This paper is full of long integrals and other technical details.
Therefore, in order to improve its readability, we sketch its conceptual
framework, which, we believe, is simple and natural. For  any DLG $(G;A,B)$  
there is a naturally defined  multiplicative unitary operator  
(in fact there are four of them). This operator is  manageable in the 
sense of Woronowicz as was proved in \cite{PS1}. Therefore, by the results 
of Woronowicz from \cite{slw}, there is a quantum group 
(modulo the problem of Haar measure) associated with any DLG.

Additionally, given a triple $(G;A,B)$, one can define two differential groupoid 
structures on $G$: $G_A$ and $G_B$ over $A$ and $B$, respectively. 
It turns out that a multiplication of one of the groupoids, say  $G_B$,  
defines a morphism in the sense of Zakrzewski \cite{SZ1} 
from $G_A$ to $G_A\times G_A$. This morphism is coassociative. Using a 
construction presented in \cite{PS2} we get a coassociative morphism from 
$C^*(G_A)$ to $C^*(G_A\times G_A)$. We may suppose that $C^*(G_A\times G_A)$ 
is (some sort of) $C^*(G_A)\mt C^*(G_A)$. Therefore, the morphism constructed as above
begins to resemble a comultiplication, one of the main ingredients  of a 
quantum group  structure. Completion of such construction of a comultiplication 
complicates the fact  that we are not able to prove the equality 
$C^*(G_A\times G_A)=C^*(G_A)\mt C^*(G_A)$. Nevertheless, we have a corresponding 
morphism of {\em reduced} $C^*$-algebras. It turns out, as one can expect,  
that this  morphism is implemented by  the multiplicative unitary. 

Due to great richness of the structure of a DLG, its other components 
enable us to find natural objects giving us almost a Hopf algebra structure. 
Moreover, they enable to construct  an invariant positive functional on 
it with a natural modular group.  To be more specific, group inverse 
implements the unitary part of an antipode; a natural cocycle on a 
groupoid implements a scaling group; a distinguished fiber gives a ``counit''; 
and Haar density on the set of units of a groupoid gives us an invariant 
functional.  We give explicit formulae for all of these mappings. 

An additional advantage of the geometric approach we use in this paper is
that, after recognizing the geometric content of a formula, we usually can prove 
it using only the very standard and well known results from calculus: the Fubini 
theorem and change of variables for a Riemann integral of smooth functions.

Although from a theoretical point of view this work adds nothing to the pure
theory of quantum groups, it shows, we believe, interesting links of the theory
of quantum groups with differential geometry, and, may be even more importantly,
provides us  with plenty of examples.  In the examples constructed this way we 
have a natural, dense $*$-subalgebra  of ``smooth''  elements in our  $C^*$-algebra, 
which, putting aside subtleties  with densities, is an algebra  of smooth, 
compactly supported  functions,  with a kind of convolution as a product.
As a result, we can use various techniques and estimates from the theory of 
integral  and appeal to our geometric intuitions in order to do our computations. 
Finally, these examples can serve as testing ground for general theory of 
non-compact  quantum groups.  Such a theory, we believe, is not finished  
even at the conceptual level. 

There are also connections of our work with ``quantization'' of  some 
semi-direct product Poisson-Lie groups, among them  $\kappa$-Poincar\'e group.  
However, in this case the situation  is not so good. Namely, the construction  
presented in the paper  does not directly apply, since the set of decomposable 
elements in $G$ i.e.   
$\{g\in G\,:\,g=ab=\ti{b}\ti{a}\,,\,a,\ti{a}\in A\,,\,b,\ti{b}\in B\}$ is only
open and dense in $G$, instead of being equal to $G$. Nevertheless, the 
methodology proposed in this paper can be applied also in this case and may 
give us correct results. 

Let us now briefly describe the content of the paper. In section 1 we collect
basic facts about differential groupoids, their morphisms and related algebras, 
in order to make the  paper more self-contained.  Section 2 is devoted to proofs 
of some  technical results used later on. In the third section we present groupoids
related to DLG. In section 4 we prove that our constructions indeed provide a 
comultiplication in the sense of the theory of quantum groups. In section 5 
we show that the algebra of a groupoid $G_A$ is very similar to a Hopf algebra, 
by identifying other ingredients of its structure and proving some needed formulae. 
The sixth section is devoted to the Haar  measure.  We use the results 
from \cite{PS2}, where it was shown that there is a natural class of KMS 
weights on any differential groupoid. Since  in a groupoid $G_A$ the set of 
units is a group $A$, one can distinguish invariant half-densities on $A$. 
These half-densities give us Haar measure on a quantum group based on $G_A$.
In the last section  we collect all the constructions and results.
Finally, there are two appendixes, where we give a geometric 
interpretation  of a cocycle implementing a scaling group, and relate an algebra
associated to a multiplicative unitary with a groupoid algebra. 

In the end of the introduction the author wants to say, that the idea
of construction presented in this paper was born during his collaboration
with Stanislaw Zakrzewski. The author also wants to thank Robert Owczarek
for his useful comments.


\section{Differential groupoids, their morphisms and algebras.}
\noindent
In this section we introduce  notation, recall basic facts about 
groupoids, their morphisms   (in the sense of S. Zakrzewski) and  
action of morphisms on  groupoid algebras. For a detailed exposition 
of the subject we refer to \cite{SZ1} (differential groupoids and 
morphisms) and \cite{PS2} (constructions related to groupoid algebras).\vs 

\noindent
{\bf The category of differential groupoids.}

Let us recall that a {\em differential relation} $r$ from a manifold $X$ 
to a manifold $Y$ is a triple $(X,Y;R)$,  such that  $R=:Gr(r)$ is a 
submanifold in $Y\times X$. By a submanifold we always mean 
{\em embedded submanifold}. All manifolds considered here are {\em smooth}, 
Hausdorff, and have  countable bases of  neighborhoods. A relation $r$ from 
$X$ to $Y$  will be denoted by $r: X \rel Y$. For a relation $r: X \rel Y$, 
by $r^T$ we denote a transposed relation, i.e. $r^T: Y\rel X$ is given by 
$Gr(r^T):=\{(x,y)\in X\times Y\,:\,(y,x)\in Gr(r)\}$.
A composition of differential relations may fail to be a differential 
relation.  Therefore, in a  definition below we assume that there are 
compositions, which  give again differential relations.  
There are also other conditions
which, together with this composition property, are contained in the notion 
of {\em transversality} \cite{SZ1}. A transversality of differential  
relations $r$ and $s$ will be denoted by $r\tran s$. A differential 
relation $r:X\rel Y$ is  {\em a differential reduction}
if $r=fi^T$, where $i$ is an inclusion mapping of a submanifold 
$C\subset X$ and  $f: C\lra Y$ is a surjective submersion. 
Now we recall a definition of a differential groupoid.
\begin{defi}{\em \cite{SZ1} Let $\G$ be a manifold.}  
A differential groupoid structure on $\G$ {\em is a triple $(m,s,e)$, 
where $m:\G\times \G\rel \G$ is a differential reduction,
$e:\{1\} \rel \G$ is a differential relation, $s:\G\lra \G$ 
is an involutive  diffeomorphism and the following 
conditions are satisfied:\vsm
\begin{enumerate}
\item $m(m\times id)=m(id\times m)$\vsm
\item  $m(e\times id)=m(id\times e)=id$\vsm
\item $m(s\times s)\sim=sm$, 
where $\sim: \G\times\G\ni(x,y)\mapsto (y,x)\in\G\times \G$\vsm
\item For any $x\in \G$  $\emptyset\neq m(x,s(x))\subset e(\{1\})$\vsm
\item $m\tran(m\times id)\,,\,m\tran(id \times m)\,,\,m\tran(e\times id)\,,\,
m\tran(id \times e)$.\vsm
\end{enumerate}}
\end{defi}
It follows that objects defined in this way coincide with standard
differential groupoids. The set $E:=e(\{1\})$ is a submanifold of $\G$ called
{\em the set of identities}.  There are also  two projections 
(surjective submersions): {\em source or right projection:} 
$\er:\G\ni x\mapsto m(s(x),x)\in E$, and {\em target or left projection:} 
$\el:\G\ni x\mapsto m(x,s(x))\in E$. A fiber of $\el$ passing through $g$ 
(i.e. the set $\{x\in\G : \el(x)=\el(g)\}$) will be denoted by
$F_l(g)$, and similarly a fiber of $\er$ by $F_r(g)$.

Now comes a definition of a morphism. It should be stressed that
this definition {\em does not coincide} with the standard one. For 
examples and motivations see \cite{SZ1,PS2}.
\begin{defi}{\em \cite{SZ1} Let $(\G,m,s,e)$ and 
$(\G',m',s',e')$ be differential groupoids. A differential relation
$h:\G\rel \G'$ is} a morphism {\em from $\G$ to $\G'$ iff\vsm
\begin{enumerate}
\item $hm=m'(h\times h)$\vsm 
\item $hs=s'h$\vsm
\item $he=e'$\vsm
\item $m'\tran(h\times h)\,,\,h\tran e$.\vsm
\end{enumerate}}
\end{defi}
In the next proposition we collect basic properties of
objects associated with a morphism.
\begin{prop}{\em \cite{SZ1}} \label{mor-prop}\notka{mor-prop}
Let $h:\G\rel\G'$ be a morphism of differential groupoids. Then
\begin{enumerate}
\item The formula $Gr(f_h):=(E\times E')\cap Gr(h^T)$ defines a smooth
mapping $f_h:E'\lra E$.
\item For each $b\in E'$, the relation $h^T$ restricted to 
$F_l(b)\times F_l(f_h(b))$ defines a smooth mapping $h_b^L:F_l(f_h(b))\lra F_l(b)$. 
The same is true for restriction to right fibers with the resulting map
$h_b^R:F_r(f_h(b))\lra F_r(b)$.
\item  The set $\G\timh\G':=\{(x,y)\in
\G\times\G'\,:\,\er(x)=f_h(\el'(y))\} $ is a submanifold of $\G\times\G'$. 
\item The set $\G*_h\G':=\{(x,y)\in \G\times \G'\,:\, \el(x)=f_h(\el'(y))\}$
is a submanifold of $\G\times\G'$.
\item The sets $\G*_h E'$ and $ \G\timh E'$ are submanifolds and the mappings
$h^L: \G*_h E'\ni (x,b)\mapsto h_b^L(x)\in \G'$ and 
$h^R: \G\timh E'\ni (x,b)\mapsto h_b^R(x)\in \G'$ are smooth.
\item The mapping $m_h: \G\timh \G'\ni(x,y)\mapsto m'(h^R(x,\el'(y)),y)\in \G'$
is  a surjective submersion.
\item The mapping $\got_h:\G\timh\G'\ni (x,y)\mapsto (x,m_h(x,y))\in
\G*_h\G'$ is a diffeomorphism.
\item The mapping: $\pi_2: \G\timh\G'\ni(x,y)\mapsto y\in \G'$ 
is a surjective submersion and $\pi_2^{-1}(y)$ is diffeomorphic to
$F_r(f_h(\el'(y)))$.
\item The mapping $\tilde{\pi}_2:\G*_h\G'\ni(x,y)\mapsto y\in\G'$ 
is a surjective submersion and $\tilde{\pi}_2^{-1}(y)$ is diffeomorphic to
$F_l(f_h(\el'(y)))$. 
\item Points 6-9 remain true if we replace $\G'$ by $F_r(y)$ for any $y\in\G'$
and restrict the corresponding mappings in an obvious way.
\end{enumerate}\vsm\dowl
\end{prop}\vsm

\noindent
{\bf Bidensities.}

Let $\G$ be a differential groupoid and let $\omh(\el),(\omh(\er))$ be a
bundle of smooth, complex-valued, half-densities along left (right) fibers of
$\G$. By $\sA(\G)$ we denote  the linear space  of compactly supported, 
smooth sections of the bundle $\,\omh(\el)\mt\omh(\er)$ \cite{Con}. 
Its elements will be called {\em bidensities} and usually denoted by $\om$. 
So $\om(x)=\lambda(x)\mt\rho(x)\in \omh T^l_x\G\mt\omh  T^r_x\G$, where
$T^l_x\G:=T_x(F_l(x))$ and $T^r_x\G:=T_x(F_r(x))$. In the following we also
write $\oml(x):=\Om^{1/2} T^l_x\G$ and $\omr(x):=\Om^{1/2} T^r_x\G$.

{\em We also use the following notation: if $M,N$ are manifolds,  
$F:M\lra N$ and $\Psi$ is some geometric object on $M$ which can be 
pushed-forward by $F$, then we denote the push-forward of $\psi$ simply 
by $F\psi$. What it really means will be clear from the context.} So below
we write for example $s(v)$ instead of $\Lambda T_xs(v)$ 
for $v\in\lma T_x^l\G$ (for a vector space $V$ by $\Lambda^{max} V$ we 
denote the maximal non-zero  exterior power of $V$).

The groupoid inverse induces a star operation on $\sA(\G)$: 
$$\om^*(x)(v\mt w):=\overline{\om(s(x))(s(w)\mt s(v))}\,,\,v\in
\Lambda^{max}T^l_x\G\,,\,w\in \Lambda^{max}T^r_x\G$$
Because $s$ is an involutive diffeomorphism which interchanges 
left and right fibers, the $*$-operation  is a well defined 
antilinear involution.

One can define  a multiplication on the vector space $\sA(\G)$. 
This multiplication introduces a $*$-algebra structure on $\sA(\G)$. 
The formula for multiplication is a special case of a more general 
construction presented in \cite{PS2} and will be given later on.
Before giving the formula, we  define some special sections of 
$\omh(\el)\mt\omh(\er)$, which are very convenient for computations. 
 
Since left(right) translations are diffeomorphisms of left(right) fibers, 
we can define in a natural way left(right)-invariant sections of 
$\omh(\el)(\omh(\er))$. Any left-invariant section of  $\omh(\el)$ is 
determined by its value on $E$ and,  conversely, any section of 
$\omh(\el)|_{E}$ can be uniquely extended  to a left-invariant section of 
$\omh(\el)$.

So let $\tilde{\lambda}$  be a non-vanishing, real, half-density on $E$ 
along left fibers (one constructs such a density by covering $E$ with maps 
adapted to the  submersion $\el$ and using appropriate partition of unity to 
glue them together). We define
$$\lambda_0(x)(v):=\tilde{\lambda}(\er(x))(s(x) v)\,,\, v\in\lma T^l_x\G.$$
This $\lambda_0$ is a left-invariant, real, non-vanishing section of $\omh(\el)$.  
Now, $\tilde{\rho}:=\tilde{\lambda} s $ is a non-vanishing, real, half-density on 
$E$ along right fibers, and $\ro$ defined by
$$\ro(x)(w):=\tilde{\rho}(\el(x)(w s(x))\,,\,w\in\lma T^r_x\G$$
is a right-invariant, non-vanishing, real section of $\omh(\el)$.
Let $\omega_0:=\lo\mt \ro$, then this is a real, non-vanishing
bidensity (of course, it doesn't belong to $\sA(\G)$ but we will still
call it bidensity).\vsm 
\begin{center}
{\em From now on the symbol $\om_0$ will always mean bidensity
constructed in this way.}\vsm
\end{center}
 
When $\om_0$ is chosen any element $\om\in\sA(\G)$ can be written uniquely
as $\omega=f\,\omega_0$ for some smooth, complex-valued function $f$ with compact
support. Note the following simple:
\begin{lem} \label{f-gw}\notka{l-gw} 
If $\om=f \om_0$ then $\om^*=f^*
\om_0$, where $f^*(x):=\overline{f(s(x))}.$
\\\dowl
\end{lem} 

\noindent
{\bf Action of groupoid morphisms on bidensities}.

It turns out that groupoid morphisms act on bidensities, i.e. 
for a morphism  $h:\G\rel\G'$ one can  construct a linear
mapping $\hat{h}:\sA(\G)\lra Lin(\sA(\G'))$ which behaves very nicely
with respect to a composition of morphisms and $*$-operation.
This construction depends in a crucial way on the prop. \ref{mor-prop}.
Let us briefly describe the construction here.

Let $a\in E'$,  $(x,y)\in\G\timh F_r(a)$\,,\,$\got_h(x,y)=:(x,z)$ 
and $b:=\el'(z)$. 
Due to the prop. \ref{mor-prop} we have the following isomorphisms:
$$i_1:\omr(x)\mt\omr(y)\lra \omh T_{(x,y)}(\G\timh F_r(a))\vsm,$$
$$\got_h:\omh T_{(x,y)}(\G\timh F_r(a))\lra \omh T_{(x,z)}(\G*_h F_r(a)),
\vsm$$ 
$$i_2:\oml(x)\mt\omr(z)\lra \omh T_{(x,z)}(\G*_h F_r(a)).$$
Thus:  $(i_2)^{-1}\got_h i_1(\rho_x\mt\rho_y)=:\lambda_x\mt\rho_z$ for some
$\lambda_x\mt\rho_z\in\oml(x)\mt\omr(z)$.\\
Moreover, since the mapping:
$F_l(y)\ni u \mapsto h_{\el'(y)}^R(x) u \in F_l(z)$
is a diffeomorphism, it defines isomorphism between $\oml(y)$ and  $\oml(z)$.\\
Now let $\om=\lambda\mt\rho\in\sA(\G)\,$  $\om'=\lambda'\mt\rho'\in\sA(\G')$. 
Then  $(i_2)^{-1} \got_h i_1(\rho(x)\mt\rho'(y))=:
\tilde{\lambda_x}\mt\tilde{\rho_z},$ 
and $h_{\el'(y)}^R(x)\lambda'(y)=:\tilde{\lambda_z'}$, so  the  expression
$[\lambda(x)\tilde{\lambda_x}]\mt\tilde{\lambda_z'} \mt\tilde{\rho_z} $
defines  a 1-density on $F_l(f_h(b))$  with values in one-dimensional 
vector space $\oml(z)\mt\omr(z)$.
In this way  we can define a mapping $\hat{h}$ mentioned above 
$$(\hat{h}(\om) \om')(z):=
\int_{F_l(f_h(b))}[\lambda\tilde{\lambda}]\mt\tilde{\lambda_z'}
\mt\tilde{\rho_z}. $$
Choose  $\om_0=\lambda_0\mt\rho_0$, $\om_0'=\lambda_0'\mt \rho_0'$. 
Then  $(i_2)^{-1} \got_h i_1 (\rho_0(x)\mt\rho_0'(y))=:
t_h(x,y)\lambda_0(x)\mt\rho_0'(z)$  for some smooth, non-vanishing function
$t_h:\G\timh\G'\lra R$ and $h_{\el'(y)}^R(x) \lambda_0'(y)=\lambda_0'(z)$.
For $\om=f_1\,\om_0\,,\,\om'=f_2\,\om_0'$  we get the explicit expression
\notka{hat-mor}
\begin{equation}\label{hat-mor}
(\hat{h}(\om)\om')(z):=\left[\int_{F_l(f_h(b))}\lo^2(x) f_1(x)
t_h(x,y)f_2(y)\right]\,\om_0'(z)=:(f_1*_hf_2)(z)\,\om_0'(z),
\end{equation} 
where $y$ is defined by $\got_h(x,y)=(x,z)$, i.e. $y=s'(h_b^L(x))z$.

Now taking $h=id:\G\rel\G$ we get algebra structure on $\sA(\G)$
and the formula for the product is \notka{product}
\begin{equation}\label{product}
 \om_1\om_2=:(f_1*f_2) \om_0\,,\,\,\,  
(f_1* f_2)(x):=\int_{F_l(x)}
\lo^2(y) f_1(y) f_2(s(y) x)=\int_{F_r(x)} \ro^2(y) f_1(x s(y)) f_2(y).
\end{equation}

It seems that there is no natural, geometric, norm on $\sA(\G)$, but one can
introduce the family of useful norms ``indexed'' by $\om_0$'s \cite{Ren,PS2}. 
Then, let us  choose $\lo$ and write $\om=f\om_0$. Let us define quantities:\notka{norm}
\begin{equation}\label{norm}
||\om||_l:=\sup_{a\in\G^0}\int_{F_l(a)}\lo^2\,|f|\,,\,\,\,
||\om||_r:=\sup_{a\in\G^0}\int_{F_r(a)}\ro^2 \,|f|\,,\,\,\, 
||\om||_0:=max\{||\om||_l,||\om||_r\}.
\end{equation}
(We do not explicitly  write the dependence on $\lo$ to 
make our  notation simpler.) It can be shown 
that $||\om||_l$, $||\om||_r$, $||\om||_0$ are norms,
and $(\sA(\G),*,||\cdot||_0)$ is a normed $*$-algebra.\vs 

\noindent
{\bf Representation of an algebra of a groupoid associated with a morphism.}

Let $\Psi$ be a smooth half-density on $\G'$ with compact support and
$\om\in\sA(\G)$, $\om=\lambda\mt\rho$.  Let $(x,y)\in \G\timh\G'$ and 
$\got_h(x,y)=:(x,z)$. As in the definition of $\hat{h}$,  
$\rho(x)\mt\Psi(y)$ can be viewed as a  half-density on 
$T_{(x,y)}(\G\timh\G')$ and  $\got_h(\rho(x)\mt\Psi(y))$ is a half-density 
on $T_{(x,z)}(\G*_h\G')$. 
Since $\oml(x)\mt\Om^{1/2}T_z\G'\simeq\Om^{1/2}T_{(x,z)}(\G*_h\G')$ 
this half-density can be written as $\tilde{\lambda_x}\mt\Psi_x(z)$ for
some $\tilde{\lambda_x}$ -- a half-density on $T_x^l(\G)$ and
$\Psi_x(z)$ -- a half-density on $T_z(\G')$. Then $\lambda(x)
\tilde{\lambda_x}\mt\Psi_x(z)$ is a 1-density on $T_x^l\G$
with values in half-densities on $T_z(\G')$. Integrating  
$\lambda(x)\tilde{\lambda_x}\mt\Psi_x(z)$ we get a half-density 
on $T_z(\G')$.  
In this way we get an operator $\pi_h(\om)$ for $\om\in\sA(\G)$ 
defined on a dense linear subspace of $L^2(\G')$
$$(\pi_h(\om)\,\Psi)(z):=\int_{F_l(f_h(b))}
[\lambda(x)\tilde{\lambda}(x)]\mt\Psi_x(z).$$
Let us choose  $\om_0$ and write $\om=f\,\om_0$. Since $\er'$ is a surjective
submersion, there is a natural isomorphism: 
$\Om^{1/2}T_w\G'\simeq\omr(w)\mt\Om^{1/2}T_{\er'(w)}E'$
for any $w\in\G'$. Therefore, if we choose $\ro'$ and
$\nu_0$ -- a non-vanishing, real half-density on $E'$, 
then $\ro'\mt\nu_0$ defines a non-vanishing, real, half-density on $\G'$. 
So any other smooth half-density with compact support $\Psi$ can be written 
as $\Psi=\psi\,\ro'\mt\nu_0=:\psi\, \Psi_0$ for
some smooth, complex-valued function $\psi$ with compact support. 
It is easy to see that: $\got_h(\ro(x)\mt\ro'(y)\mt\nu_0(a))=t_h(x,y)
(\lo(x)\mt\ro'(z)\mt\nu_0(a)),$ where $t_h$ is as in the definition of $\hat{h}$.
So the  explicit formula is\notka{rep-mor}
\begin{equation}\label{rep-mor}
(\pi_h(\om)\Psi)(z)=
\left[\int_{F_l(f_h(b)}\lo^2(x)\,f(x)\,t_h(x,y)\,\psi(y)\right]\Psi_0(z),
\end{equation}
where $b:=\el'(z)\,,\,\got_h(x,y)=(x,z)$. 
In the next proposition we collect essential properties of the mappings
$\hat{h}$ and $\pi_h$.
\begin{prop}{\em \cite{PS2}}\label{r-m-prop}\notka{r-m-prop}
Let $\G,\G',\G''$ be differential groupoids, $h:\G\rel\G'$ and $k:\G'\rel\G''$ 
differential groupoid morphisms.
For any $\om\in\sA(\G),\om',\om_1'\in\sA(\G'),\om''\in\sA(\G'')$:
\begin{itemize}
\item[a)] $\hat{k}((\hat{h}(\om)\om'))\om''=\hat{kh}(\om)(\hat{k}(\om')\om'')$
\item[b)] $\pi_k(\hat{h}(\om)\om')=\pi_{kh}(\om)\pi_k(\om')$
\item[c)] $(\om')^*(\hat{h}(\om)\om_1')=(\hat{h}(\om^*)\om')^*\om_1'$
\item[d)] $\pi_h$ is a $*$-representation of $\sA(\G)$ by bounded operators
on $L^2(\G')$, moreover for any $\om_0$ we have the estimate:
$||\pi_h(\om)||\leq||\om||_0$ (on the right hand side is used the norm defined 
in eq. \ref{norm}).
\item [e)] Let us choose some $*$-invariant $\om_0'$ and let $||\cdot||_0'$ 
denote the associated norm. For any $\om'\in\sA(\G')$ there exists a sequence
$\om_n\in\sA(\G)$ such that 
$\lim_{n\rightarrow\infty}||\hat{h}(\om_n)\om'-\om'||_0'=0$.
\end{itemize}\vsm
\dowl
\end{prop}
Due to these properties one can define a $C^*$-norm on $\sA(\G)$, then
complete $\sA(\G)$ in this norm and get a kind of  universal $C^*$-algebra 
of a differential groupoid.  
For any morphism $h$, a  mapping $\hat{h}$ can be 
uniquely extended to a morphism between corresponding 
$C^*$-algebras and $\pi_h$ to a  non-degenerate representation \cite{PS2}. 

It is also easy to see  that the representation $\pi_{id}$ is faithful on
$\sA(\G)$, so the  function $\om\mapsto ||\pi_{id}(\om)||$ defines 
a $C^*$-norm on $\sA(\G)$. The completion of $\sA(\G)$ in this 
norm is called the {\em reduced} $C^*$-algebra of $\G$ and will be
denoted by $C^*_r(\G)$.\vs

\noindent
{\bf Bisections and their action on groupoid algebras.}

Recall that a submanifold $B\subset \G$ is {\em a bisection} iff
$\el|_B$ and $\er|_B$ are diffeomorphisms. A set of bisections 
of a groupoid $\G$ is a group under a natural multiplication of
subsets of $\G$ and inverse given by $B^{-1}:=s(B).$ Bisections act
(from the left) on  $\G$  by diffeomorphisms: for $g\in\G$ $B g=g' g$, where
$g'$ is a unique point in $B$ composable with $g$. One can immediately
verify that $B F_r(g)=F_r(g)$ and $B F_l(g)=F_l(B g).$ These facts enable
us to define an action of bisections on $\sA(\G)$ according to the formula
$$(B\om)(B g)(B v\mt B w):=\om(g)(v\mt w)\,,\,
{\rm where}\,\,v\in\lma T^l_g\G\,,\,
w\in \lma T^r_g\G.$$
It is easy to see that for $\om_0=\lo\mt\ro$ and $B\om_0=: f\om_0$, 
the function $f$ is given by\vsm \notka{fun-B}
\begin{equation}\label{fun-B}f(B g) \ro(B g)(B w)=\ro(g)(w)\,,\,
w\in \lma T^r_g\G.
\end{equation}
It turns out that if $B$ is a bisection of $\G$ and $h:\G\rel \G'$ is a 
morphism, then the set $h(B)$ is a bisection of $\G'$. In the next lemma
we collect basic properties of action of morphisms on bisections which
allow us to interpret bisections as multipliers on $C^*(\G)$ and $C^*_r(\G)$.
\begin{prop}{\rm \cite{PS2}}\label{bis-prop}\notka{bis-prop}
Let $\G,\G'$ be differential groupoids, $B$ a bisection of $\G$ and 
$h:\G\rel \G'$ a morphism. 
For any $\om,\om_1\in \sA(\G)$ and $\om'\in\sA(\G')$:
\begin{enumerate}
\item $\om^*(B \om_1)=(s(B)\om)^* \om_1$\vsm
\item $(\hat{h}(B\om))\om'=h(B)(\hat{h}(\om)\om')$\vsm
\item $\pi_h(B\om)=h(B)\pi_h(\om)$.\vsm
\end{enumerate}
\dowl 
\end{prop}


\section{$C^*$-algebra of Cartesian product of groupoids.}

In this section we collect some technical results about a $C^*$-algebra of
the Cartesian product of differential groupoids. We are not able to prove that
$C^*(\G_1\times\G_2)$ is equal to (some kind of) $C^*(\G_1)\mt C^*(\G_2)$. 
However, we have the following (expected)  result.
\begin{prop}
Let $\G_1,\G_2$ be differential groupoids. Then 
$C^*_r(\G_1\times \G_2)=C^*_r(\G_1)\mt_{\sigma} C^*_r(\G_2)$, where 
$\mt_{\sigma}$ denotes the minimal tensor product.
\end{prop}
{\em Proof:}
We adopt the following notation: $\sA_i:=\sA(\G_i)$, $H_i:=L^2(\G_i)$, 
$\pi_i:\sA_i\lra B(H_i)$ -- the identity representation of $\sA_i$, 
$||\cdot||_i$ -- an operator norm on $B(H_i)$,
$A_i:=C^*_r(\G_i)=\overline{\pi_i(\sA_i)}^{||\cdot||_i}$, 
$i=1,2$ and the  corresponding objects for $\G_1\times\G_2$: 
$\sA:=\sA(\G_1\times \G_2)$, 
$H:=L^2(\G_1\times\G_2)=H_1\mt H_2$,
$\pi: \sA\lra B(H)$, $||\cdot||$ -- operator norm on  $B(H)$ and 
$A:=C^*_r(\G_1\times \G_2)=\overline{\pi(\sA)}^{||\cdot||}$. 
It is clear that $\sA_1\mt\sA_2\subset\sA$ and
$\pi|_{\sA_1\mt\sA_2}=\pi_1\mt\pi_2$.

Let us  choose  $\lo$ and $\tilde{\lo}$ -- non-vanishing, real, left-invariant 
half-densities along left fibers  on $\G_1$ and $\G_2$ respectively. 
Then $\lo\mt\ti{\lo}$ is  a non-vanishing, real, left-invariant half-density along left 
fibers on $\G_1\times\G_2$,  and we have the corresponding right-invariant  
half-densities and $*$-invariant bidensities: 
$\om^1_0=\lo\mt\ro,\om^2_0=\ti{\lo}\mt\ti{\ro}$ and $\om_0=\om^1_0\mt\om^2_0$.
Let $||\cdot||_0$ be the corresponding norm on $\sA$. We begin by the following
\begin{lem} $\sA_1\mt\sA_2$ is dense in $\sA$ in the topology defined by 
$||\cdot||_0$.
\end{lem}
{\em Proof:} Let $\sD(\cdot)$ denotes compactly supported, smooth functions. 
The following result is standard.\\
For any  $f\in \sD(\GG)$ there exist compact sets 
$K_1\subset \G_1\,,\,K_2\subset \G_2$ such that\\
$\forall \, \epsilon >0 \,\exists \,g\in\sD(\G_1)\mt\sD(\G_2)\,:\, 
{\rm supp} (g)\subset K_1\times K_2\,\,{\rm and}\,\, \sup|f-g|<\epsilon$.

Let $\om\in\sA(\GG)$, then $\om=f\om_0$ for some 
$f\in\sD(\GG)$. Let $K_1,K_2$ be subsets as above and choose some 
$h\in\sD(\GG)$, such that $h(x)\geq 0$, $h=1$ 
on $K_1\times K_2$ and let  $M:=||h\om_0||_0.$  Take $\epsilon>0$ and let
$g=\sum f_i\mt k_i\in\sD(\G_1)\mt\sD(\G_2)$ be such that 
$\sup|f-g|<\epsilon/M$. Define $\om_i^1:=f_i\om_0^1\,,\,\om_i^2:=k_i\om_0^2$.
Then simple computations show that 
$||\om-\sum \om_i^1\mt\om_i^2||_0\leq \epsilon.$\\
\dowl

Since for $\om\in \sA$ there is an inequality  $||\om||\leq ||\om||_0$, 
we get that the closure of  $\pi(\sA_1\mt\sA_2)$ contains the closure of 
$\pi(\sA)$ i.e.  $A\subset \overline{\pi(\sA_1\mt\sA_2)}^{||\cdot||}=
\overline{\pi_1(\sA_1)\mt\pi_2(\sA_2)}^{||\cdot||}$. This means that 
$A\subset A_1\mt_{\sigma} A_2$.
On the other hand for  $A_1\ni a_1=\lim \om_n,\,\om_n\in\sA_1$ and 
$A_2\ni a_2=\lim \om_n',\, \om_n'\in\sA_2$ we immediately verify that
$\pi_1(a_1)\mt\pi_2(a_2)=\lim\pi_1(\om_2)\mt\pi_2(\om_n')=
\lim \pi(\om_n\mt\om_n')$. In this way 
$\pi_1(A_1)\mt\pi_2(A_2)\subset \overline{\pi(\sA)}=C^*_r(\G_1\times\G_2)$
and $A_1\mt_{\sigma} A_2=\overline{\pi_1(A_1)\mt\pi_2(A_2)}
\subset C^*_r(\G_1\times\G_2)$.\\
\dow

\noindent
From now on, if  $A$ and  $B$ are $C^*$-algebras, a  tensor product $A\mt B$ 
means the {\em minimal tensor product}.\\
In the next lemma $\om(\om_1\mt I)$ and $\om(I\mt\om_1)$ mean product in
$M(C^*_r(\G\times \G))$ (i.e. multiplier algebra).
\begin{lem}\label{mnoz-jedno}\notka{mnoz-jedno}
Let $\om=F(\om_0\mt\om_0)\in\sA(\G\times \G)$ and $\om_1=f\om_0\in\sA(\G)$.
Then   $\om(\om_1\mt I)$ and $\om(I\mt\om_1)$ are elements of 
$\sA(\G\times \G)$, which are given by the formulae
$$\om(\om_1\mt I)=:(F*(f\mt I))(\om_0\mt\om_0)\,\,,\,\,\,
\om(I\mt\om_1)=:(F*(I\mt f))(\om_0\mt\om_0)$$ 
$$(F*(f\mt I))(g_1,g_2):=\int_{F_l(g_1)}\lo^2(g) F(g,g_2)f(s(g) g_1),$$
$$(F*( I\mt f))(g_1,g_2):=\int_{F_l(g_2)}\lo^2(g) F(g_1,g)f(s(g) g_2).$$
\end{lem}
{\em Proof:} To verify these equations it is enough to show that: 
$(F*(f\mt I))*(f_1\mt f_2)=F*(f*f_1\mt f_2)$ and 
$(F*(I\mt f))*(f_1\mt f_2)=F*(f_1\mt f* f_2)$.
Let us compute the LHS of the first equality:
$$(F*(f\mt I))*(f_1\mt f_2)(g_1,g_2)=\int_{F_l(g_1)\times F_l(g_2)}
(\lo^2(g)\mt\lo^2(h))(F*(f\mt I))(g,h)f_1(s(g) g_1)f_2(s(h) g_2)=$$
$$=\int_{F_l(g_1)\times F_l(g_2)}(\lo^2(g)\mt\lo^2(h))
f_1(s(g) g_1)f_2(s(h) g_2)
\int_{F_l(g)}\lo^2(k) F(k,h)f(s(k) g).$$
Since $F_l(g)=F_l(g_1)$ the above integral  is equal to
$$\int_{F_l(g_1)\times F_l(g_1)\times F_l(g_2)}(\lo^2(g)\mt\lo^2(k)\mt\lo^2(h))
f_1(s(g)  g_1)f_2(s(h) g_2)F(k,h)f(s(k) g).$$
And the RHS:
$$F*(f*f_1\mt f_2)(g_1,g_2)=\int_{F_l(g_1)\times F_l(g_2)}
(\lo^2(g)\mt\lo^2(h))F(g,h)(f*f_1)(s(g)  g_1)f_2(s(h) g_2)=$$
$$=\int_{F_l(g_1)\times F_l(g_2)}(\lo^2(g)\mt\lo^2(h))F(g,h)f_2(s(h) g_2)
\int_{F_r(g_1)}\ro^2(k)f(s(g) g_1 s(k)) f_1(k)=$$
$$=\int_{F_l(g_1)\times F_l(g_2)}(\lo^2(g)\mt\lo^2(h))F(g,h)f_2(s(h) g_2)
\int_{F_l(e_R(g_1))}\lo^2(k)f(s(g) g_1\gdot k)) f_1(s(k))=$$
$$=\int_{F_l(g_1)\times F_l(g_2)}(\lo^2(g)\mt\lo^2(h))F(g,h)f_2(s(h)  g_2)
\int_{F_l(g_1)}\lo^2(k)f(s(g)  k)) f_1(s(k)  g_1)=$$
$$=\int_{F_l(g_1)\times F_l(g_1)\times 
F_l(g_2)}(\lo^2(k)\mt\lo^2(g)\mt\lo^2(h))
f_1(s(k) g_1)f_2(s(h)  g_2)F(g,h)f(s(g)  k);$$
the passage from the  second to the third line follows from the fact that
$s$ restricted to $F_r(g_1)$ is a diffeomorphism onto $F_l(\er(g_1))$ and
$s(\ro)=\lo$; the equality of the third and the fourth line
is implied by diffeomorphism $F_l(\er(g_1))\ni k \mapsto g_1 k\in F_l(g_1)$.

The equality $(F*(I\mt f))*(f_1\mt f_2)=F*(f_1\mt f* f_2)$ can be proved in
the same way.\\
\dowl\\
Later on we will need the following technical result
\begin{lem}\label{diff-W}\notka{diff-W}
Let $\Phi:\G\times \G\rightarrow \G\times \G$ be a diffeomorphism,
and $\psi:\G\times \G\rightarrow R$ be a smooth, non-vanishing function.
Let $\lo$ be a  left-invariant, real, non-vanishing half-density 
and $\om_0$ corresponding *-invariant
bidensity. Then $\lo\mt\lo$ is a left-invariant, real, non-vanishing half-density on
$\G\times \G$ and $\om_0\mt\om_0$ is the corresponding *-invariant bidensity.
Let $||\cdot||_l,\,||\cdot||_r,\,||\cdot||_0$ be the norms on 
$\sA(\G\times \G)$ associated with $\lo\mt\lo$ as defined in eq. \ref{norm}.
Define $S:\sA(\G)\mt\sA(\G)\lra \sA(\G\times\G)$ by a formula:
$$ S\left({\scriptstyle \sum} \om_i\mt\om_i'\right):=
\psi\, \Phi\left({\scriptstyle \sum} f_i\mt f_i'\right) 
(\om_0\mt\om_0),$$
where $\om_i=:f_i\om_0\,,\,\om_i'=:f_i'\om_0$ and $(\Phi g)(x):=g(\Phi^{-1}(x))$ 
is  a push-forward of a function $g$.
Then $S(\sA(\G)\mt\sA(\G))$ is dense in $\sA(\G\times\G)$ in the topology 
defined by $||\cdot||_0$.
\end{lem}
{\em Proof:} Let $\om=:F (\om_0\mt\om_0)\in\sA(\G\times\G)$ and 
$K:={\rm supp\,} F$. 
Define $\ti{F}:=F/\psi$ and $F_1:=\Phi^{-1}\ti{F}$. Then 
${\rm supp\,} F_1=\Phi^{-1}(K)$ is a compact set. There exist compact sets 
$K_1,K_2$ such that  $\Phi^{-1}(K)\subset K_1\times K_2$ and 
$$\forall \epsilon \,\exists f_i,g_i : {\rm supp\,} f_i\subset K_1\,,\,
{\rm supp\,} g_i\subset K_2\,,\,\sup|F_1-{\scriptstyle \sum}f_i\mt g_i|\leq \epsilon$$
Then ${\rm supp\,} \Phi({\scriptstyle \sum}f_i\mt g_i)
\subset \Phi(K_1\times K_2)$.
Choose a smooth function $h\in \sD(\G\times \G)$  such that 
$0\leq h\leq 1$ and $h=1$ on $K\cup \Phi(K_1\times K_2)$.
Let $M:=||h||_0$ and $N:=
\sup \{|\psi(g_1,g_2)|:(g_1,g_2)\in K\cup \Phi(K_1\times K_2)\}$.
Let $\epsilon>0$ be given and let $f_i,g_i$ be as above, and 
$\sup|F_1-{\scriptstyle \sum}f_i\mt g_i|\leq \epsilon/{M N}$.
Now for $\ti{\om}:={\scriptstyle \sum}f_i\om_0\mt g_i\om_0$ we have:
$$||\om-S(\ti{\om})||_l=\sup_{(e_1,e_2)\in E\times E}\int_{F_l(e_1,e_2)}
\lo^2(g_1)\mt\lo^2(g_2)
|F(g_1,g_2)-\psi(g_1,g_2) \Phi({\scriptstyle \sum}f_i\mt g_i)(g_1,g_2)|.$$
We can estimate the integral as follows:
$$\int_{F_l(e_1,e_2)}\lo^2(g_1)\mt\lo^2(g_2)
|F(g_1,g_2)-\psi(g_1,g_2) \Phi({\scriptstyle \sum}f_i\mt g_i)(g_1,g_2)|=$$
$$=\int_{F_l(e_1,e_2)}\lo^2(g_1)\mt\lo^2(g_2)
|h(g_1,g_2)||F(g_1,g_2)-\psi(g_1,g_2) \Phi({\scriptstyle \sum}f_i\mt g_i)(g_1,g_2)|=$$
$$=\int_{F_l(e_1,e_2)}\lo^2(g_1)\mt\lo^2(g_2)|h(g_1,g_2)|
|\psi(g_1,g_2)||\ti{F}(g_1,g_2)-\Phi({\scriptstyle \sum}f_i\mt g_i)(g_1,g_2)|\leq$$
$$\leq N \sup |\ti{F}-\Phi({\scriptstyle \sum}f_i\mt g_i)|
\int_{F_l(e_1,e_2)}\lo^2(g_1)\mt\lo^2(g_2)
|h(g_1,g_2)|.$$
Since 
$$\sup |\ti{F}-\Phi({\scriptstyle \sum}f_i\mt g_i)|=
\sup |(\Phi F_1)-\Phi({\scriptstyle \sum}f_i\mt g_i)|=
 \sup |F_1-{\scriptstyle \sum}f_i\mt g_i)|\leq\epsilon/{M N}.$$
we get an  inequality:$||\om-S(\ti{\om})||_l\leq \epsilon$.
In the same way we can estimate $||\om-S(\ti{\om})||_r$.\\
\dowl


\section {Double Lie Groups and related objects.}

In this section we introduce basic objects of our investgations --
groupoids related to a Double Lie Group (DLG).
We begin by recalling the definition of a DLG
(also known as a matched pair of Lie groups or 
a bicrossproduct Lie group).
\begin{defi} {\em \cite{WeLu}} A double Lie group {\em is a triple of 
Lie groups $(G;A,B)$  such that $A,B$ are closed subgroups of 
$G$, $A\cap B=\{e\}$ and $G=AB$.}
\end{defi}
The structure of a DLG  defines four projections: 
$$\al,\ar : G\lra A\,,\,\,  \bl,\br :G\lra B\,\,{\rm by:}\,\, 
g=\al(g)\br(g)=\bl(g)\ar(g).$$ We also define a  relation 
$m_A: G\times G\rel G$ by:
$$Gr(m_A):=\{(b_1 a b_2; b_1 a,a b_2)\,:\,a\in A\,,\,b_1,b_2\in B\},$$ 
and a diffeomorphism  
$s_A:G\ni g\mapsto \ar(g)\br(g)^{-1}=\bl(g)^{-1}\al(g)\in G$.\\
It turns out \cite{SZ1} that $G_A:=(G,m_A,A,s_A)$ is a differential groupoid.
The same is true for $G_B:=(G,m_B,B,s_B)$, where $m_B,s_B$ are defined similarly
to $m_A,s_A$. Moreover $\d:=m_B^T :G_A\rel G_A\times G_A$ and 
$s_B:G_A\rel G_A$ are  morphisms of differential groupoids. In the ``extreme case''
$A=\{e\},B=G$ we have $G_A=G$ (the  groupoid $G_A$  is a group) and $G_B=(G,diag^T,G,id)$,
where $diag:G\ni g\lra (g,g)\in G\times G$ is the diagonal mapping. The basic
example of a DLG, the reader may think of, is an Iwasawa decomposition of 
a semi-simple Lie group $G=K(AN)$.

In the next lemma we explicitly describe mappings
and sets related to the morphism $\d$. We use the same notation as in 
prop. \ref{mor-prop}. The proof is straightforward.
\begin{lem}\label{roznewzorki}\notka{roznewzorki} 
Let $(G;A,B)$ be a double Lie group and let $\d:=m_B^T$. Then: 
\begin{enumerate}
\item $Gr(\d)=\{(a_2 b_2, b_2 a_3;a_2 b_2 a_3): a_2,a_3 \in A,b_2\in B\};$
\item $f_{\d}:A\times A\ni(a_1,a_2)\mapsto a_1a_2\in A;$
\item $G_A\times_{\d}(G_A\times G_A)=\{(\ti{b_1} a_2 a_3;a_2 b_2,a_3 b_3) : 
\ti{b_1},b_2,b_3\in B,a_2,a_3\in A\};$
\item $G_A\times_{\d}F_r(\ti{a_2},\ti{a_3})=
\{(\ti{b_1} \al(\ti{b_2}\ti{a_2})\al(\ti{b_3}\ti{a_3});\ti{b_2}\ti{a_2},
\ti{b_2}\ti{a_2}) :  \ti{b_1},\ti{b_2},\ti{b_3}\in B\};$
\item $G_A*_{\d}(G_A\times G_A)=
\{(a_2 a_3 b_1 ; a_2 b_2,a_3 b_3) : b_1,b_2,b_3\in B,a_2,a_3\in A\};$
\item $G_A *_{\d}F_r(\ti{a_2},\ti{a_3})=
\{(\al(\ti{b_2}\ti{a_2})\al(\ti{b_3}\ti{a_3}) b_1;\ti{b_2}\ti{a_2},
\ti{b_2}\ti{a_2}) : b_1,\ti{b_2},\ti{b_3}\in B\};$
\item $\d^L(g_1;a_2,a_3)=(a_2\bl(a_3\br(g_1)),a_3 \br(g_1))\,\,, \al(g_1)=a_2 a_3;$
\item $\d^R(g_1;a_2,a_3)=(\bl(g_1)a_2,\br(\bl(g_1)a_2) a_3)\,\,, \ar(g_1)=a_2 a_3;$
\item $m_{\d}(g_1;g_2,g_3)=(\bl(g_1)g_2, \br(\bl(g_1)\al(g_2)) g_3);$
\item $\got_{\d}(g_1;g_2,g_3)=(g_1;\bl(g_1)g_2, \br(\bl(g_1)\al(g_2)) g_3).$
\end{enumerate}\dowl
\end{lem}

For any DLG one can define a pentagonal diffeomorphism 
$W:G\times G\lra G\times G$:\notka{W}
\begin{equation}\label{W}
W(s,t):=(s\al(t)^{-1},\br(s\al(t)^{-1})t)\,\,,\,\,\,\,
W^{-1}(s,t)=(s \al(\br(s)^{-1} t), \br(s)^{-1} t).
\end{equation}
By a push-forward of half-densities, $W$  defines a multiplicative unitary
operator on $L^2(G)\mt L^2(G)$ , which will also be denoted by $W$.
Then $W^*$ is a push-forward by $W^{-1}$.

Now we give an interpretation of $W$ in terms of groupoids $G_A$ and $G_B$.
We need this  to prove  easily 
that $W\in M(CB(L^2(G))\mt C^*_r(G_A))$ ($CB(H)$ stands for compact operators).
However, this interpretation 
can be also used to show that $W$ is a unitary bicharacter on 
quantum groups with some universal properties and to construct a 
quantum double. 
For a groupoid $\G=(\G,m,E,s)$
we denote by $\G^{op}$ the groupoid with the reverse multiplication i.e.
$(\G,m^{op},E,s)$. As it was shown in \cite{SZ2}, the pentagonal diffeomorphism
$W$ is equal to $(id\times m_A)(m_B^T\times id)$ and the pentagon equation
follows from the fact that $m_B^T$ is a morphism $G_A\rel G_A\times G_A$.
(Now it is clear that  there are four multiplicative operators, since
one can interchange $A$ and $B$; one can also consider $G_A^{op}$ and $G_B^{op}$.)
But $W$ is not only a diffeomorphism but  a diffeomorphism of a very
special kind, namely, it is implemented by a bisection.
Consider the set $U:=\{(g,s_B(g))\,:\,g\in G\}$. This is
clearly a submanifold in $G\times G$. The following lemma is a result of 
direct computations.
\begin{lem}\label{W-bis}\notka{W-bis}
$U$ is a bisection of $G_B^{op}\times G_A$, moreover $U\cdot(s,t)=W(s,t)$.
(the left hand side of this equality is understood as an action of a bisection
on a groupoid element).\\
\dowl
\end{lem}
Because  a bisection of a differential groupoid $\Gamma$ is an element
of $M(C^*_r(\Gamma))$ we conclude that 
$W\in M(C^*_r(G_B^{op}\times G_A))=M(C^*_r(G_B^{op})\mt C^*_r( G_A))$. 
The algebra  $C^*_r(G_B^{op})$ acts on $L^2(G_B^{op})=L^2(G)$ in a 
non-degenerate way, therefore $W\in M(CB(L^2(G))\mt C^*_r(G_A))$.\vs

In what follows we  will also need various modular functions related to 
DLG, so now we fix notation:\notka{fun-mod}
\begin{equation}\label{fun-mod}
\Adaa(a)=\det (P_A Ad(a) P_A)\,\,,\,\,\,
\Adab(a)=\det (P_B Ad(a) P_B)\,,\,a\in A;\end{equation}\vsm\vsm\vsm
\begin{equation}\Adba(b)=\det (P_A Ad(b) P_A)\,,\,\,
\Adbb(b)=\det (P_B Ad(b) P_B)\,,\,\,b\in B,
\end{equation}
where  \prb\,and \pra\,  denote projections in $\g$ 
corresponding to the  decomposition $\g=\ab\ms\bb$.

The last object we define  is a smooth function $Q: G\lra R\setminus\{0\}$.\\
Let $Ad(g)=:\left( \begin{array}{cc} g_1 & g_2\\
g_3 & g_4 \end{array} \right)$ be a decomposition of operators of 
the adjoint representation with respect to the direct sum structure
$\g=\ab\ms\bb$ i.e.
$g_1:\ab\lra\ab\,,\,g_2:\bb\lra\ab$, etc.
Then $Q$ is defined by\notka{defQ}
\begin{equation}\label{defQ}
G\ni g\mapsto Q(g):=
\frac{\det (Ad(g))}{\det(g_1)\det(g_4)}\in R\setminus\{0\}.
\end{equation}
For $b\in B,a\in A$ we have: 
 $Ad(a)=:\left(\ba{cc} \alpha_1 & \alpha_2\\0 & \alpha_4\ea\right)
\,,\,Ad(b)=:\left(\ba{cc} \beta_1 & 0\\\beta_3 & \beta_4\ea\right)$, and
$$ Q(ba)=\frac{\det(\beta_4 \alpha_4)}
{\det(\beta_3\alpha_2+\beta_4\alpha_4)}. \label{qba}$$
The function $Q$ is related to the modular functions by equalities:
\notka{Q-mod}
\begin{equation}\label{Q-mod}
Q(g)=\frac{\Adaa(\al(g))\Adba(\br(g))}{\Adaa(\ar(g))\Adba(\bl(g))}=
\frac{\Adab(\ar(g))\Adbb(\bl(g))}{\Adab(\al(g))\Adbb(\br(g))}.
\end{equation}
The relationship of the function $Q$ with groupoids $G_A$ and $G_B$ is 
described in the following lemma, the proof of which is straightforward.
\begin{lem}\label{koc}\notka{koc} $Q$ is a one cocycle on $G_A$ and $G_B$, i.e.
$Q(ab)Q(ba')=Q(aba')\,,\,Q(ba)Q(a b')=Q(bab')$ for any 
$a,a'\in A\,,\,b,b'\in B$. Moreover, $Q$ is invariant with respect to the group
inverse i.e. $Q(g)=Q(g^{-1})$.
\dowl
\end{lem}
The function $Q$ is exactly the one which appears in the definition of 
manageability of $W$ \cite{slw, PS1}. 
For the geometric meaning of $Q$ see Appendix.\vs

\noindent
{\bf Multiplication in $\sA(G_A)$.}

Let us choose some  $\mu_0\neq 0$ -- a real half-density on $T_eB$ 
and define\notka{lambda}
\begin{equation}\label{lambda}
\lo(g)(v):=\mu_0(g^{-1}v)\,,\,\,v\in \lma T^l_gG_A.
\end{equation}
It is easy to see that this is a left-invariant, non-vanishing, half-density 
on $G_A$ and the  corresponding right-invariant half-density is given by\notka{ro}
\begin{equation}\label{ro}
\ro(g)(w)=|\Adab(\al(g))|^{-1/2}\mu_0(wg^{-1})\,,\,w\in \lma T^r_gG_A
\end{equation} 
We put $\om_0=\lo\mt\ro$, and from the formula (\ref{product}) we obtain
the expression for multiplication in $\sA(G_A)$.
$$((f_1\om_0)(f_2\om_0))(g)=:(f_1*f_2)(g)\om_0(g)$$\vsm\vsm\vsm\vsm
\notka{product-Ga}
\begin{eqnarray}\label{product-Ga}
(f_1*f_2)(g)=\int_B \mu^2_L(b) \,f_1(\al(g)b)\,f_2(b_L(\al(g)b)^{-1}g)=\\
=\int_B \mu^2_R(b) \,|\Adab(\al(ba_R(g)))|^{-1}\,f_1(gb_R(ba_R(g))^{-1})\,
f_2(ba_R(g)),
\end{eqnarray} 
where  $\mu_L$ and $\mu_R$ are  left and right-invariant  half-densities  
on $B$  defined by $\mu_0$ (we use the fact that left and right fibers are 
diffeomorphic to $B$).\\
We finish this section with a simple observation, which will be  used
later on.
\begin{lem}\label{mb-bis}\notka{mb-bis}
For $g_0\in G$ the sets $\{g\in G : (g_0,g)\in \d(G)\}=\br(g_0)A$ and 
$\{g\in G : (g,g_0)\in \d(G)\}=A\bl(g_0)$ are  bisections of 
$G_A$.\\
\dowl
\end{lem}



\section{Comultiplication}

This section is entirely devoted to a proof that indeed $\hat{\d}$ can be 
used to obtain a comultiplication in the sense of the theory of quantum groups.
\begin{tw}\label{main} \notka{main} 
Let $(G; A, B)$ be a double Lie group,  $G_A=(G,m_A,s_A,A)\,,\,
G_B=(G,m_B,s_B,B)$ corresponding differential groupoids and let $\d:=m_B^T$. Then
\begin{itemize}
\item[a)] The mapping $\hat{\d}$ extends uniquely to 
$\Delta\in Mor(C^*_r(G_A), C^*_r(G_A)\mt C^*_r(G_A))$;
\item[b)] $\Delta$ is coassociative, i.e. 
$(\Delta\mt id) \Delta=(id \mt \Delta)\Delta$;
\item[c)] For any $a,b\in C^*_r(G_A)$,  elements    
        $\Delta(a)(b\mt I)$ and $\Delta(a)(I\mt b)$ belong 
$C^*_r(G_A)\mt C^*_r(G_A)$. Moreover, the linear spaces  
        $span\{\Delta(a)(b\mt I): a,b\in C^*_r(G_A)\}$ and 
$span\{\Delta(a)( I\mt b): a,b\in C^*_r(G_A)\}$ are dense in 
$C^*_r(G_A)\mt C^*_r(G_A)$.
\end{itemize}
\end{tw}

In general we know \cite{PS2} that  $\hat{\d}$ extends to a morphism from
$C^*(G_A)$ to $C^*(G_A\times G_A)$ and it is not a priori clear that it
defines also a morphism between reduced algebras. However, it is easy to see
that, for $\om\in\sA(G_A)$,  $\hat{\d}(\om)$ defines a multiplier on
$C^*_r(G_A\times G_A)$. 

Indeed, let $\pi_{id}$ denotes the identity representation
of $\sA(G_A\times G_A)$ on $L^2(G_A\times G_A)$. Using the prop. 
\ref{r-m-prop} we have the inequality 
$||\pi_{id}(\hat{\d}(\om)\om_1)||=||\pi_{\d}(\om)\pi_{id}(\om_1)||
\leq ||\pi_{\d}(\om)||||\pi_{id}(\om_1)||$. From this inequality we infer 
that  $\hat{\d}(\om)$ is a  bounded, linear mapping, defined 
on a dense, linear subspace of $C^*_r(G_A\times G_A)$. 
Therefore, it can be extended  in a unique way to the whole 
$C^*_r(G_A\times G_A)$. Moreover, since
$\om_2^*(\hat{\d}(\om)\om_1)=(\hat{\d}(\om^*)\om_2)^*\om_1$, we see that 
$\hat{\d}(\om)$ is adjointable. Thanks to this property,
 $\hat{\d}(\om)$ defines a multiplier. 

Let us choose some $\om_0'\in\sA(G_A\times G_A)$ and let
$\om\in\sA(G_A\times G_A)$. Take  $\om_n\in\sA(G_A)$  as in  item e) 
of the prop.\ref{r-m-prop}. We have an estimate:
$||\pi_{id}(\hat{\d}(\om_n)\om)-\pi_{id}(\om)||
\leq ||\hat{\d}(\om_n)\om-\om||_0'$ and, since this tends to $0$, we see
that any element of $\sA(G_A\times G_A)$ can be approximated in the
norm defined by identity representation by elements from
$\hat{\d}(\sA(G_A))\sA(G_A\times G_A)$. So the same is true for the whole
$C^*_r(G_A\times G_A)$.

In this way to prove statement a) of the theorem  we need
continuity of $\hat{\d}$ as a mapping defined on a dense subspace
$\sA(G_A)\subset C^*_r(G_A)$. This will immediately follow from the
following.
\begin{prop}\label{implement}\notka{implement}
Let  $W$ be the  pentagonal diffeomorphism defined in (\ref{W}).   
Representation $\pi_\d$ is implemented by $W$ i.e.
$$\pi_\d(\om)=W(\pi_{id}(\om)\mt I)W^*$$
\end{prop}
{\em Proof:} Of course this result should be expected taking
into account the close relationship between the  algebra defined by $W$ in
a ``standard way'' and $C^*_r(G_A)$ (see Appendix B). 
We start by giving a formula for $\hat{\d}$. For $\om=:f\om_0\in\sA(G_A)$ 
and $\om_1=:F(\om_0\mt\om_0)\in \sA(G_A\times G_A)$, 
by (\ref{hat-mor}) and lemma \ref{roznewzorki}  we have 
$\hat{\d}(\om)\om_1=(f*_{\d}F)(\om_0\mt\om_0)$, and
$$(f*_\d F)(a_1b_1,a_2b_2)=\int_{F_l(a_1a_2)} \lo^2(a_1a_2b)
f(a_1 a_2 b)\, 
t_\delta(a_1a_2b;\bl^{-1}(a_1a_2 b) a_1 b_1,\ar(a_2 b) b^{-1} b_2)\,\times$$
$$\times F(\bl^{-1}(a_1a_2 b) a_1 b_1,\ar(a_2 b) b^{-1} b_2)$$
As it was proved in \cite{PS2}, for any morphism $h$ of differential
groupoids,  the function $t_h$ is right-invariant with respect
to multiplication in the second groupoid i.e. $t_h(x,y)=t_h(x,\el(y))$.
In the next lemma we prove a formula for this function.
\begin{lem} $t_\d(b_1 a_2 a_3;a_2,a_3)=|\Adbb(\br(b_1a_2))|^{-1/2}$.
\end{lem}
{\em Proof: } We use notation and results given in prop \ref{mor-prop} and
lemma \ref{roznewzorki}. 
Let $x:=b_1a_2 a_3\,,\, y:=(a_2,a_3)$ and $z:=m_\d(x,y)=(b_1a_2,\br(b_1a_2)a_3)$.
Let  $W:=T_{(x,y)}(G_A\times_\d F_r(y))$  and  $U:=\ker \pi_2\subset W$.
Vectors from $U$  are represented by curves 
$u(t):=(b_1(t)a_2a_3;a_2,a_3)\,,\, b_1(0)=b_1$. 
It is easy to see that a  subspace $V\subset W$ of 
vectors represented by curves 
$$v(t):=(b_1\al(b_2(t)a_2)\al(b_3(t)a_3);
b_2(t)a_2, b_3(t)a_3)\,,\,b_2(0)=b_3(0)=e$$
is complementary to $U$.
So the isomorphism $i_1:\omr(x)\mt\omr(y)\lra \omh W$ is given by \notka{i1}
\begin{equation}\label{i1}
i_1(\ro(x)\mt\ro'(y))(u\dz v):=\ro(x)(\pi_1u)\ro'(y)(\pi_2v)\,,\,
u\in \lma U,v\in\lma V.
\end{equation}
Let $\ti{W}:=T_{(x,z)}(G_A*_\d F_r(y))$.
Then $\ti{U}:=\ker\ti{\pi}_2\subset \ti{W}$ is a subspace of vectors 
represented by curves 
$\ti{u}(t):=(b_1a_2a_3b_1(t);  b_1a_2,\br(b_1a_2)a_3)\,,\, b_1(0)=e$. 
Choosing some $\ti{V}\subset\ti{W}$ complementary to $\ti{U}$ we can
write isomorphism $i_2:\oml(x)\mt\omr(z)\lra \omh \ti{W}$ as \notka{i2}
\begin{equation}\label{i2}
i_2(\lo(x)\mt\ro'(z))(\ti{u}\dz\ti{v}):=
\lo(x)(\ti{\pi}_1\ti{u})\ro'(z)(\ti{\pi}_2\ti{v})\,,\ti{u}\in\lma\ti{U},
\ti{v}\in\lma\ti{V}.
\end{equation}
The function $t_\d$ is defined by $\got_\d i_1(\ro(x)\mt\ro'(y))=:
t_\d(x,y) i_2(\lo(x)\mt\ro'(z))$, i.e.\notka{i12}
\begin{equation}\label{i12}
i_1(\ro(x)\mt\ro'(y))(u\dz v)=
t_\d(x,y)i_2(\lo(x)\mt\ro'(z))(\got_\d(u)\dz \got_\d(v)).
\end{equation}
Since the subspace $\got_\d(V)$ is represented by curves
$$t\mapsto (b_1\al(b_2(t)a_2)\al(b_3(t)a_3);
b_1b_2(t)a_2, \br(b_1\al(b_2(t)a_2))b_3(t)a_3),$$
one immediately sees  that it is  complementary to $\ti{U}$. 
So to compute $t_\d$ we  need  a  decomposition of $\got_\d(u)$
with respect to the direct sum $\ti{W}=\ti{U}\ms \got_\d(V)$. 
This  is given by curves $\ti{u}(t), v(t)$ such that 
$\got_\d u(t)=(\got_d v(t))(x^{-1},z^{-1})\ti{u}(t)$ (group multiplication in
$G\times G\times G$). Using our parametrization we see that for a given curve 
$b_1(t)\,,\,b_1(0)=b_1$ we have to find curves 
$\ti{b}_1(t)\,,\,\ti{b}_1(0)=e$  and
$b_2(t),b_3(t)\,,\,b_2(0)=b_3(0)=e$ such that
$$(b_1(t)a_2a_3; b_1(t)a_2,\br(b_1(t)a_2)a_3)=$$
$$=(b_1\al(b_2(t)a_2)\al(b_3(t)a_3);b_1b_2(t)a_2,
\br(b_1\al(b_2(t)a_2))b_3(t)a_3)$$
$$(b_1a_2a_3;b_1a_2,\br(b_1a_2)a_3)^{-1}
\,(b_1a_2a_3\ti{b}_1(t);b_1a_2,\br(b_1a_2)a_3)$$
The solution of this  equation is given by \notka{rozw}
\begin{equation}\label{rozw}
\ti{b}_1(t)=b_R(b_1^{-1}b_1(t)a_2a_3)\,,\,
b_2(t)=b_1^{-1}b_1(t)\,,\,
b_3(t)=b_R(b_1^{-1}b_1(t)a_2).
\end{equation}

Now we know enough  to perform computations for 
the function $t_\d$.
Let $(X_1,\dots,X_n)$ be a basis in $T_eB$ such that 
$\mu_0(X_1\wedge\dots\wedge X_n)=1$. Then $(X_1x,\dots,X_n x)$ is a basis
in $T^r_xG_A$ and $u_i:=(X_i x;0_{a_2},0_{a_3})\,,i=1,\dots,n$ 
form a basis in $U$. We have also corresponding
bases in $T^r_{a_2}G_A$ and  $T^r_{a_3}G_A$, and 
$(v_1,\dots,v_n,v_1',\dots,v_n')$, where 
$v_i:=(b_1\al(X_i a_2)a_3;X_i a_2,0_{a_3})$ and
$v_i':=(b_1a_2\al(X_i a_3);0_{a_2}, X_i a_3)$, is a basis in $V$. 
For $\ro'=\ro\mt\ro$, $\ro$ as in (\ref{ro}), the equation
(\ref{i1}) gives
$$i_1(\ro(x)\mt\ro'(y))(u_1\wedge\dots\wedge v_n')=
|\Adab(\al(x))\Adab(a_2)\Adab(a_3)|^{-1/2}.$$
Let $P$ be a projection onto $\ti{U}$ corresponding to the decomposition
$\ti{W}=\ti{U}\ms \got_\d(V)$. Using (\ref{i2}) we can now rewrite (\ref{i12})
as 
$$|\Adab(\al(x))\Adab(a_2)\Adab(a_3)|^{-1/2}=
t_\d(x;a_2,a_3) i_2(\lo(x)\mt\ro'(z))(\got_\d(u_1)
\wedge\dots\wedge\got_\d(u_n)\wedge\dots \wedge\got_\d(v_n'))=$$
$$=t_\d(x;a_2,a_3)
\lo(x)(\ti{\pi}_1P\got_\d u_1\wedge\dots\wedge\ti{\pi}_1P\got_d u_1)
\ro'(z)(\ti{\pi}_2\got_\d v_1\wedge\dots\wedge\ti{\pi}_2\got_\d v_n').$$ 
From (\ref{rozw}) we get   
$$\ti{\pi}_1Pu_i=x \br(x^{-1}X_ix)\,\,,\,\,
\ti{\pi}_2\got_\d v_i=(b_1X_ia_2,\br(b_1\al(X_ia_2))a_3)\,\,,\,\,
\ti{\pi}_2\got_\d v_i'=(0_{b_1a_2},\br(b_1a_2)
X_ia_3),$$ 
therefore  
$$\lo(x)(\ti{\pi}_1P\got_\d u_1\wedge\dots\wedge\ti{\pi}_1P\got_d u_1)\,
\ro'(z)(\ti{\pi}_2\got_\d v_1\wedge\dots\wedge\ti{\pi}_2\got_\d v_n')=$$
$$=|\Adab(\al(b_1a_2))
\Adab(\al(\br(b_1a_2)\ti{a_3}))|^{-1/2}\times$$
$$\times|\det (P_BAd(x^{-1})P_B)|^{1/2}|\det (P_BAd(b_1)P_B) |^{1/2}
|\det (P_BAd(\br(b_1a_2))P_B) |^{1/2}.$$
But  since $\Adab(a_1a_2)=\Adab(a_1)\Adab(a_2)$ and 
$\al(b_1a_2)\al(\br(b_1a_2)a_3)=\al(b_1a_2a_3)=\al(x)$,
our equation for $t_\d$ reduces to
$$|\Adab(a_2 a_3)|^{-1/2}=t_\d(x;a_2,a_3)
|\det (P_BAd(x^{-1})P_B Ad(b_1)P_B Ad(\br(b_1a_2))P_B)|^{1/2}.$$
Since $x=b_1a_2 a_3$ we have
$$P_BAd(x^{-1})P_B Ad(b_1)P_B Ad(\br(b_1a_2))P_B=
P_BAd(a_3^{-1}a_2^{-1})P_B Ad(\br(b_1a_2))P_B,$$
and finally we get $t_\d(b_1a_2a_3;a_2,a_3)=|\Adbb(\br(b_1a_2))|^{-1/2}$.\\
\dowl\vs

\noindent
In this way we arrive at the formula for $\hat{\d}$: (we recall
that $f*_\d F$ is defined by 
$\hat{\d}(f\om_0)(F(\om_0\mt\om_0))=:(f*_\d F)(\om_0\mt\om_0)$) \notka{komno}
\begin{equation}\label{komno}
%
(f*_{\d}F)(a_1b_1,a_2 b_2)=$$
$$=\int_{F_l(a_1a_2)} \lo^2(a_1a_2b) |\Adbb(\bl(a_2 b))|^{-1/2}
f(a_1 a_2 b)\, F(\bl^{-1}(a_1a_2 b) a_1 b_1,\ar(a_2 b) b^{-1} b_2)=$$
$$=\int_{B} \mu_L^2(b) |\Adbb(\bl(a_2 b))|^{-1/2}
f(a_1 a_2 b)\, F(\bl^{-1}(a_1a_2 b) a_1 b_1,\ar(a_2 b) b^{-1} b_2).
\end{equation}

Now we can prove the main formula of the prop.\ref{implement}.
Let $\ro$ be the right-invariant half-density defined in 
(\ref{ro}). Choose $\nu_0$ -- a real half-density on $T_eA$ and let $\nu_l$ be the 
corresponding left-invariant half-density on $A$. Then $\psi_0:=\ro\mt\nu_l$ is
a non-vanishing half-density on $G$.  Define $\Psi_0:=\psi_0\mt\psi_0$, 
then this is a non-vanishing half-density on $G\times G$. 
The explicit formula for $\Psi_0$ is
$$\Psi_0(s,t)(Ys\dz sX\dz Y_1t\dz t X_1)=$$
$$=\ro(s)(Ys)\ro(t)(Y_1t)
\nu_l(\ar(s))(\ar(s)X)\nu_l(\ar(t))(\ar(t)X_1)=$$
$$=|\Adab(\al(s)\al(t))|^{-1/2}\mu_0(Y)\mu_0(Y_1)\nu_0(X)\nu_0(X_1)\,,
\,\,Y,Y_1\in \lma T_eB\,,\,X,X_1\in \lma T_eA.$$
In the next lemma we compute an action of $W$ on $\Psi_0$.
\begin{lem}\notka{W'}\label{W'} 
Let $(s,t)\in G\times G$,  $(\ti{s},\ti{t}):=W^{-1}(s,t)$, 
$X,X_1\in \lmA$ and  $Y,Y_1\in\lmB$.
$$\Adbb(\br(s))\Adaa(\al(\ti{t}))\, 
W^{-1}(Y s \wedge s X \wedge Y_1 t \wedge t X_1)= 
(Y \ti{s} \wedge \ti{s} X \wedge Y_1 \ti{t} \wedge \ti{t} X_1).$$
\end{lem}
{\em Proof: }
For $(s,t)\in G\times G$ let  
$i_{st}:\bb \ms  \ab \ms \bb \ms  \ab\lra T_{(s,t)}(G\times G)$ be
an isomorphism given by the decomposition 
$T_{(s,t)}(G\times G)=\bb s\ms s \ab \ms \bb t\ms t \ab$. Then  
$$W^{-1}(Y s \wedge s X \wedge Y_1 t \wedge t X_1)=
\det(i_{\ti{s}\ti{t}}^{-1} T_{(s,t)}(W^{-1})i_{st}) 
(Y \ti{s} \wedge \ti{s} X \wedge Y_1 \ti{t} \wedge \ti{t} X_1).$$
Let us denote the composition 
$i_{\ti{s}\ti{t}}^{-1} T_{(s,t)}(W^{-1})i_{st}$ by $F$ and 
write $F$ in a block form corresponding to the decomposition 
$T_{(e,e)}(G\times G)=\bb \ms  \ab \ms \bb \ms  \ab$:
$$F=\left(\begin{array}{cccc}F_{11} & F_{12} & F_{13} & F_{14}\\
                            F_{21} & F_{22} & F_{23} & F_{24}\\
                            F_{31} & F_{32} & F_{33} & F_{34}\\         
                            F_{41} & F_{42} & F_{43} & F_{44}
        \end{array}\right).$$
From the formula for $W^{-1}$ given in eq. (\ref{W}) we easily obtain:
$$F_{11}=I=F_{44}\,\,\,{\rm  and }\,\,\,
F_{12}=F_{13}=F_{14}=F_{41}=F_{42}=F_{43}=0,$$
so  $\det F=\det  \left(\begin{array}{cc} F_{22} & F_{23}\\
                                          F_{32} & F_{33}\\             
        \end{array}\right)$, 
and  mappings $F_{22}, F_{23},F_{32}, F_{33}$ are given by:
$$F_{22}: \ab\ni X\mapsto 
Ad(\al(z)^{-1})X-\pra Ad(\br(z)\al(t)^{-1})\prb Ad(\br(s)^{-1})X,$$
$$F_{23}: \bb\ni Y\mapsto \pra Ad(\br(z)\al(t)^{-1})Y,$$
$$F_{32}: \ab\ni X\mapsto -Ad(\br(s)^{-1})\prb Ad(\br(s)^{-1})X,$$
$$F_{33}: \bb\ni Y\mapsto Ad(\br(s)^{-1})Y,$$
where $z:=\br(s)^{-1}\al(t)$. 

\noindent
Now it remains to compute the determinant and compare it with the definitions 
of $\Adbb$ and $\Adaa$.\\
\dowl\\
Having these results we are ready to prove prop. \ref{implement}.
Let a function $w:G\times G\lra R$ be defined by 
$W\Psi_0=:w \Psi_0$ and  $W(F)$ by: $W(F\Psi_0)=:W(F)\Psi_0$.
Then $W(F)(s,t)=F(W^{-1}(s,t)) w(s,t)$.
From the preceding lemma we get
$$w(s,t)=|\Adbb(\br(s)) \Adaa(\al(\ti{t})) \Adab(\al(\ti{t}))|^{-1/2}
\,\,{\rm where}\,\,\ti{t}=\br(s)^{-1}t.$$
Now let us compute $\pi_{\d}(\om)W$. For $\om=:f\om_0$ and $\Psi=:F\Psi_0$
using the formula (\ref{komno}) we get
$$(f*_{\d}W(F))(a_1 b_1,a_2 b_2)=
\int_B \mu_L^2(b) |\Adbb(\bl(a_2 b))|^{-1/2}\,f(a_1 a_2 b)\, 
W(F)(\bl^{-1}(a_1 a_2 b) a_1 b_1, \ar(a_2 b) b^{-1} b_2)=$$
$$=\int_B \mu_L^2(b) |\Adbb(\bl(a_2 b))|^{-1/2}\,f(a_1 a_2 b)
w(\bl^{-1}(a_1 a_2 b) a_1 b_1, \ar(a_2 b) b^{-1} b_2)\times$$ 
$$\times F(\bl^{-1}(a_1 a_2 b) a_1 a_2 \bl(a_2^{-1} b_1),b_1^{-1} a_2 b_2)=$$
$$=|\Adaa(\al(b_1^{-1} a_2)) \Adab(\al(b_1^{-1} a_2))|^{-1/2}\times$$
$$\times\int_B \mu_L^2(b) |\Adbb(\bl(a_2 b)) \Adbb(\bl^{-1}(a_2 b) b_1)|^{-1/2}
f(a_1 a_2 b) F(\bl^{-1}(a_1 a_2 b) a_1 a_2 \bl(a_2^{-1} b_1) , 
b_1^{-1} a_2 b_2)=$$ 
$$=|\Adaa(\al(b_1^{-1} a_2)) \Adab(\al(b_1^{-1} a_2)) \Adbb(b_1)|^{-1/2}
\int_B \mu_L^2(b)f(a_1 a_2 b)\, 
F(\bl^{-1}(a_1 a_2 b) a_1 a_2 \bl(a_2^{-1} b_1) ,  b_1^{-1} a_2 b_2)$$
On  the other, hand using the formula for multiplication (\ref{product-Ga}) and
lemma \ref{mnoz-jedno} we can write
$$W((f\mt I)*F)(a_1 b_1 a_2 b_2)=w(a_1 b_1, a_2 b_2) ((f\mt I)*F)
(a_1 b_1 \al(b_1^{-1} a_2), b_1^{-1} a_2 b_2)=$$
$$=|\Adaa(\al(b_1^{-1} a_2)) \Adab(\al(b_1^{-1} a_2)) \Adbb(b_1)|^{-1/2}
\int_B \mu_L^2(b) f(a_1 a_2 b)\, 
F(\bl^{-1}(a_1 a_2 b) a_1 a_2 \bl(a_2^{-1} b_1) ,  b_1^{-1} a_2 b_2)$$
Comparing both expression we see that
$\pi_{\d}(\om)W=W(\pi_{id}(\om)\mt I)$. The proposition \ref{implement} is 
proved.\\\dow\\

In this way the  proof of the statement a) of the theorem \ref{main} 
is completed i.e. $\hat{\d}$ can be extended to a morphism $\Delta$
from $C^*_r(G_A)$ to $C^*_r(G_A)\mt C^*_r(G_A)$.\vs\\
b) The prop. \ref{implement} immediately yields coassociativity 
of $\Delta$ by a  standard proof based on the pentagonal equation.\vs\\
c) Density conditions.

Our proof of the density conditions will be the following. 
First, let us note that, by using standard
density arguments, one easily shows that if 
$\hat{\d}(\om)(\om_1\mt I)\in \sA(G_A\times G_A)$ then also for any
$a,b\in C^*_r(G_A)$ we have $\Delta(a)(b\mt I)\in C^*_r(G_A)\mt C^*_r(G_A)$
and the same holds for the second inclusion.

We begin by  proposing the explicit formulae for
$\hat{\d}(\om)(I\mt \om_2)$ and $\hat{\d}(\om)(\om_1\mt I)$ for
$\om,\om_1,\om_2\in\sA(G_A)$ based on some geometric considerations.
These objects will be elements of $\sA(G_A\times G_A)$. 
Next, we  verify that our guess is  correct. Finally,
we check the density conditions. 

Let us  first give  some heuristic arguments for our formulae. 
Let us write $\hat{\d}(\om)=\sum\om_k'\mt\om_k''$ and  formally compute
$$\hat{\d}(\om)(I\mt\om_1)(g_1,g_2)=\sum\om_k'(g_1)(\om_k''\om_1)(g_2)=
\left[\left(\sum\om_k'\mt\om_k''\right)(g_1,\cdot)\om_1\right](g_2)=
[\hat{\d}(\om)(g_1,\cdot)\om_1](g_2)$$
What kind of object is $\hat{\d}(\om)(g_1,\cdot)$? 
Since $\d=m_B^T$ then we can  expect that 
$\hat{\d}(\om)(g_1,\cdot)=\om(m_B(g_1,\cdot))$. Now using lemma 
\ref{mb-bis}, we can interpret $\om(m_B(g_1,\cdot))$ {\em as a function
on a bisection} $\br(g_1)A$ and this can be identified with a multiplier 
on $\sA(G_A)$. Therefore, the natural candidate for $\hat{\d}(\om)(I\mt\om_1)$ is 
$ \om(m_B(g_1,\cdot)) \om_1$. It follows that (almost) the same is true
in the  differential case.
\begin{re}{\em The author knew the formula for $\hat{\d}$ (\ref{komno}) before
the formula for $\hat{h}$ for an arbitrary morphism of differential
groupoids given in section 1 was introduced in \cite{PS2}.
The formula for $\hat{\d}(\om)(I\mt\om_1)$ was found as below and then
the equality $\hat{\d}(\om)(\om_1\mt\om_2)=
(\hat{\d}(\om)(I\mt\om_2))(\om_1\mt I)$ was used. }
\end{re}
Let $\d^R$ and $\d^L$ be the mappings associated to the morphism $\d$ 
described in lemma \ref{roznewzorki}:\\
$\d^L(g_1;a_2,a_3)=(a_2\bl(a_3\br(g_1)),a_3 \br(g_1))\,\,, \al(g_1)=a_2 a_3$,
\\
$\d^R(g_1;a_2,a_3)=(\bl(g_1)a_2,\br(\bl(g_1)a_2) a_3)\,\,, \ar(g_1)=a_2 a_3$.\\
In the following lemma we denote the tangents to these mappings by the 
same letters.
\begin{lem}\label{dl-dr}\notka{dl-dr}
For any $(g_1,g_2;g)\in \d$ the mappings:
$\d^L:T^l_gG_A\lra T^l_{(g_1,g_2)}(G_A\times G_A)=
T^l_{g_1}G_A\ms T^l_{g_2}G_A$ and 
$\d^R:T^r_gG_A\lra T^r_{(g_1,g_2)}(G_A\times G_A)=T^r_{g_1}G_A\ms T^r_{g_2}G_A$ 
are injective.\\
Moreover, the mappings:
$\pi_1\d^L:T^l_gG_A\lra T^l_{g_1}G_A\,\,,\,\,
\pi_1\d^R:T^r_gG_A\lra T^r_{g_1}G_A\,\,,\,\,
\pi_2\d^L:T^l_gG_A\lra T^l_{g_2}G_A\,\,,\,\,
\pi_2\d^R:T^r_gG_A\lra T^r_{g_2}G_A$,
where $\pi_1,\pi_2$ are projections from $G_A\times G_A$ to $G_A$, are isomorphisms.
\end{lem}
{\em Proof:} Let $g_1:=a_1 b_1=\ti{b}_1\ti{a}_1\,,\,g_2:= b_1 \ti{a}_2=a_2 b_2\,,\,
g:=\d^T(g_1,g_2)=a_1 b_1 \ti{a}_2=\ti{b}_1\ti{a}_1\ti{a}_2=a_1 a_2 b_2$. 
Let $X\in T^l_gG_A$ and $Y\in T^r_gG_A$ be represented by curves 
$a_1 a_2 b_2(t)\,,\,b_2(0)=b_2\,,\,\ti{b}_1(t)\ti{a}_1\ti{a}_2\,,\,
\ti{b}_1(0)=\ti{b}_1$ respectively.
Then $\d^L(X)$ and $\d^R(Y)$ are represented by $(a_1\bl(a_2b_2(t)),a_2b_2(t))$ and 
$(\ti{b}_1(t)\ti{a}_1,\br(\ti{b}_1(t)\ti{a}_1)\ti{a}_2)$ 
and the first statement is clear.
 Direct computations also prove that 
$(\pi_1\d^L)^{-1}$, $(\pi_2\d^L)^{-1}$, $(\pi_1\d^R)^{-1}$ and 
$(\pi_2\d^R)^{-1}$ are given by: \\
$(\pi_1\d^L)^{-1}: a_1 b_1(t)\mapsto a_1 a_2 \bl(a_2^{-1}b_1(t))\,,\,b_1(0)=b_1\,$;
$\,\,(\pi_2\d^L)^{-1}: a_2 b_2(t)\mapsto a_1 a_2 b_2(t)\,,\,b_2(0)=b_2$\\
$(\pi_1\d^R)^{-1}: \ti{b_1}(t)\ti{a_1}\mapsto 
\ti{b_1}(t)\ti{a_1}\ti{a_2}\,,\,\ti{b}_1(0)=\ti{b}_1\,$;
$\,\,(\pi_2\d^R)^{-1}: \ti{b_2}(t)\ti{a_2}\mapsto \br(\ti{b_2}(t)\ti{a_1}^{-1})
\ti{a_1}\ti{a_2}\,,\,\ti{b}_2(0)=b_1$
\\\dowl\\
Now suppose we are given $\om\in\sA(G_A)$,  $g_1\in G_A$ and
$u_1:=v_1\mt w_1\in \lma T^l_{g_1}G_A\mt \lma T^r_{g_1}G_A.$
For each $g_2$ such that $(g_1,g_2)\in\d(G_A)$ we can define a number:
$$\d_1(\om;g_1,u_1)(g_2):=
\om(\d^T(g_1,g_2))((\pi_1\d^L)^{-1}v_1\mt (\pi_1\d^R)^{-1}w_1).$$
Since the set of possible $g_2$ is the  bisection  $b_R(g_1)A$, 
in this way we get a  function on this bisection.
In the same way  for each $g_2\in G_A$ and
$ u_2:=v_2\mt w_2\in \lma T^l_{g_2}G_A\mt \lma T^r_{g_2}G_A.$
we get a  function on the bisection $A b_L(g_2)$ given by
$$\d_2(\om;g_2,u_2)(g_1):=
\om(\d^T(g_1,g_2))((\pi_2\d^L)^{-1}v_2\mt (\pi_2\d^R)^{-1}w_2).$$
Having a bisection $U$ and a function $f$ on it, we can define a natural 
action of  this pair on bidensities by  the formula  
$((f U)\om)(g):=f(g') (U\om)(g)$, where $g'$ is the unique point in $U$ 
such that $\al(g')=\al(g)$.
Now for $\om,\om_1,\om_2\in \sA(G_A)$ we define
$K_1(\om,\om_1)\in\sA(G_A\times G_A)$ and $K_2(\om,\om_2)\in\sA(G_A\times G_A)$ by
\begin{eqnarray}
K_1(\om,\om_1)(g_1,g_2)(u_1\mt u_2):=
[\d_1(\om,g_1,u_1) (\br(g_1)A) \om_1](g_2)(u_2)\\
K_2(\om,\om_2)(g_1,g_2)(u_1\mt u_2):=
[\d_2(\om,g_2,u_2) (A\bl(g_2))  \om_2](g_1)(u_1).
\end{eqnarray}
It is clear that these formulae are bilinear, so in fact 
$K_1,K_2:\sA(G_A)\mt\sA(G_A)\lra \sA(G_A\times G_A)$. 
These mappings are our candidates for $\hat{\d}(\om)(I\mt\om_1)$ and
$\hat{\d}(\om)(\om_1\mt I)$.
In the next lemma we compute  functions $\d_1$ and $\d_2$.
\begin{lem}\label{d1d2} \notka{d1d2}
Let $(g_1,g_2)\in\d^T(G_A)\,$,$u_1,u_2$ be as above. 
Then the functions $\d_1$ and $\d_2$ are given by:
$$\d_1(f \om_0; g_1, u_1)(g_2)=
f(\d^T(g_1,g_2))|\Adab(\al(g_2))|^{-1/2}|\Adab(\ar(g_2))|^{-1/2}
\om_0(g_1)(u_1).$$
$$\d_2(f \om_0; g_2, u_2)(g_1)=
f(\d^T(g_1,g_2))\om_0(g_2)(u_2).$$
\end{lem}
{\em Proof:} Let $g_1=a_1 b_1\,,\,g_2=b_1 \tilde{a}_2\,,\,
g:=\d^T(g_1,g_2)=a_1 b_1 \tilde{a}_2$ and $u_1=v_1\mt w_1$. 
Using the definitions of $\om_0$ and $\d_1$ we compute
$$\d_1(\om;g_1,u_1)(g_2):=
\om(g)((\pi_1\d^L)^{-1}v_1\mt (\pi_1\d^R)^{-1}w_1)
=f(g) 
\om_0((\pi_1\d^L)^{-1}v_1\mt (\pi_1\d^R)^{-1}w_1)=$$
$$=f(g)\mu_0(g^{-1}(\pi_1\d^L)^{-1}v_1)|\Adab(\al(g))|^{-1/2}
\mu_0((\pi_1\d^R)^{-1}w_1 g^{-1}).$$
Let us represent elements from $T^l_{g_1}G_A$ by $g_1X$ for $X\in T_eB$ and
elements from $T^r_{g_1}G_A$ by $Xg_1$ for $X\in T_eB$. Then from the 
lemma \ref{dl-dr} we see that the mapping  
$T_eB\ni X\mapsto g^{-1}(\pi_1\d^L)^{-1}g_1 X\in T_eB$ is equal to 
$P_BAd(\ti{a}_2^{-1})$, and the mapping  
$T_eB\ni X\mapsto (\pi_1\d^R)^{-1}Xg_1 g^{-1}\in T_eB$ is equal to the 
identity mapping. Since $\al(g)=\al(g_1) \al(g_2)$ the first formula is proved.
In the same way the second one can be proved.\\\dowl\vs

\noindent
To find explicit expressions for  $K_1$ and $K_2$ we also  need   
formulae for actions  of bisections $bA$ and $Ab$ on $\sA(G_A)$.
\begin{lem}\label{bis-b-lp}\notka{bis-b-lp} For  $b_0\in B$ actions of 
bisections $b_0A$ and $Ab_0$ on $\om_0$ are given by:
$$((b_0A)\om_0)(g)=|\Adbb(b_0)|^{-1/2} \left|\frac{\Adab(\al(b_0^{-1}g))}
{\Adab(\al(g))}\right|^{-1/2} \om_0(g)\,\,\,{\rm and}\,\,\, 
((A b_0)\om_0)(g)=|\Adbb(b_0)|^{-1/2}\om_0(g).$$
\end{lem}
{\em Proof:}
Let $(b_0A)\om_0=:f \om_0$, $h:=(b_0A)\gdot g=b_0 g$. 
Using the equation  (\ref{fun-B}) we get
$$f(h)=\left|\frac{\Adab(\al(g))}{\Adab(\al(h))}\right|^{-1/2}
\frac{\mu_0(w)}{\mu_0(((b_0A)\gdot (w g))h^{-1})}\,\,
{\rm for }\,\,w\in\lma T_eB.$$
Let $X\in T_eB$ be represented by a curve $b(t)$, 
then $((b_0A)\gdot (X g))h^{-1}$ is represented
by $b_0b(t)b_0^{-1}$ and this gives us the first formula.\vs\\
To prove the second one, let us again put  $h:=(Ab_0) \gdot g$ and 
$(Ab_0)\om_0=:f\om_0$, but now
$$h=b_L(\al(g)b_0^{-1})^{-1} g=a_R(\al(g)b_0^{-1})b_0 b_R(g)\,\,{\rm  and}
\,\,g=b_L(\al(h) b_0)^{-1}h=a_R(\al(h)b_0)b_0^{-1}b_R(h).$$
For $X$ represented by a curve $b(t)$, the vector  
$((Ab_0)\gdot (X g))h^{-1}$ is  represented by
$$\ti{b}(t):=\bl(\al(b(t)g)b_0^{-1})^{-1}b(t)\bl(\al(g)b_0^{-1})$$
To simplify expressions let us write:  $g=:a b$. Now we observe that  
$\ti{b}(t)=:b_1(t)b_2(t)$, for 
$$b_1(t):=\bl(\al(b(t)a)b_0^{-1})^{-1}\bl(ab_0^{-1})\,,\,b_1(0)=e\,\,
{\rm  and }\,\,  
b_2(t):=\bl(ab_0^{-1})^{-1} b(t)\bl(ab_0^{-1})\,,\,b_2(0)=e.$$ 
One easily verifies that 
$\,\,b_1(t)=
\br[\br(b_0a^{-1})a\ar(a^{-1}b(t)^{-1}a)a^{-1}\br(b_0a^{-1})^{-1}]$.\newline
From this decomposition it follows that the tangent mapping is given by 
$$T_eB\ni X\mapsto Ad(\br(b_0a^{-1}))X-P_BAd(\br(b_0a^{-1})a)P_AAd(a^{-1})X$$
Let 
$Ad(a)=:\left(\begin{array}{cc}\alpha_1 & \alpha_2\\0&\alpha_4\end{array}
\right)$ 
and $Ad(\br(b_0a^{-1}))=:
\left(\begin{array}{cc}\beta_1 & 0\\\beta_3&\beta_4\end{array}\right)$ be the 
decompositions as in (\ref{defQ}).
With this notation, our mapping is equal to 
$\beta_4+\beta_3\alpha_2\alpha_4^{-1}$.
Substituting it to the expression for $f$ and using the formula for $h$
one gets
$$f(h)=\left|\frac{\Adab(\al(g))}{\Adab(\ar(ab_0^{-1}))}\right|^{-1/2}
|\det(\beta_4+\beta_3\alpha_2\alpha_4^{-1})|^{-1/2}=
\left|\frac{\det(\beta_4\alpha_4 +\beta_3\alpha_2)}
{\Adab(\ar(ab_0^{-1}))}\right|^{-1/2}$$
But since $\br(b_0a^{-1})a=\ar(ab_0^{-1})b_0$, we get 
$f(h)=|\Adbb(b_0)|^{-1/2}$.\\\dowl

Now, using  lemmas \ref{bis-b-lp} and \ref{d1d2}, we  easily derive
the formulae for the mappings $K_1$ and $K_2$  
(we use notation 
$K_i(f \om_0,f_i\om_0)=:K_i(f,f_i)(\om_0\mt\om_0)$).\notka{K1-K2}
\begin{eqnarray}\label{K1-K2}
K_1(f,f_1)(a_1 b_1,a_2 b_2)= |\Adab(\al(b_1^{-1} a_2))|^{-1}|\Adbb(b_1)|^{-1/2}
f(a_1 b_1\al(b_1^{-1} a_2)) f_1(b_1^{-1} a_2 b_2)\\
\label{K1-K2-2}
K_2(f,f_2)(a_1 b_1, a_2 b_2)=|\Adbb(\bl(g_2))|^{-1/2} 
f(a_1 a_2 b_2) f_2(\bl^{-1}(a_1\bl(g_2)) a_1 b_1).
\end{eqnarray}
\begin{lem}\label{K-del}\notka{K-del}
$K_1(\om,\om_1)(\om_2\mt I)=\hat{\d}(\om)(\om_2\mt\om_1)=K_2(\om,\om_2)(I\mt\om_1)$.
\end{lem}
{\em Proof:} Using the formula (\ref{K1-K2}) 
and lemma \ref{mnoz-jedno} we compute 
$$K_1(f,f_1)*(f_2\mt I)(a_1 b_1,a_2 b_2)=$$
$$=\int_B\mu_L^2(b) |\Adbb(b)|^{-1/2} |\Adab(\al(b^{-1}a_2))|^{-1}
f(a_1 b \al(b^{-1}a_2))  f_2(\bl^{-1}(a_1 b) a_1 b_1) f_1(b^{-1} a_2 b_2)$$
The mapping $B\ni b\mapsto \bl(a_2^{-1}b)\in B$ is a diffeomorphism. 
Applying it to the integrated density one gets exactly eq. (\ref{komno}).
The second equality is immediate from lemma \ref{mnoz-jedno}.\\
\dowl\vs\\
Now we know enough to get rid of $K_1$ and $K_2$ i.e. identify them
with $\hat{\d}(\om)(I\mt\om_1)$ and $\hat{\d}(\om)(\om_1\mt I)$.
\begin{prop}\label{out-K}\notka{out-K} 
Let  $\om=f\om_0,\om_1=f_1\om_0\in \sA(G_A)$. Then 
$\hat{\d}(\om)(I\mt\om_1)\,,\,\hat{\d}(\om)(\om_1\mt I)\in\sA(G_A\times G_A)$.
Moreover, for $f*_{\d}(I\mt f_1)$ and $f*_{\d}(f_1\mt I)$ given by:
$$\hat{\d}(f\om_0)(I\mt f_1\om_0)=:(f*_{\d}(I\mt f_1)) (\om_0\mt\om_0)\,\,
{\rm  and}\,\, 
\hat{\d}(f\om_0)(f_1\om_0\mt I)=:(f*_{\d}(f_1\mt I))(\om_0\mt\om_0)$$  we have
the following expressions:
$$(f*_{\d}(I\mt f_1))(g_1,g_2)= |\Adab(\al(\br(g_1)^{-1}g_2))|^{-1}
|\Adbb(\br(g_1))|^{-1/2}
f(g_1\al(\br(g_1)^{-1} g_2)) f_1(\br(g_1)^{-1}g_2)$$
$$(f*_{\d}(f_1\mt I))(g_1,g_2)=|\Adbb(\bl(g_2))|^{-1/2} 
f(\al(g_1) g_2) f_1(\bl(\al(g_1)g_2)^{-1} g_1).$$
\end{prop}
{\em Proof:} This proposition simply states that 
$\hat{\d}(\om)(I\mt\om_1)=K_1(\om,\om_1)$ and 
$\hat{\d}(\om)(\om_1\mt I)=K_2(\om,\om_1)$. Let us verify the first equality.
We know  that $K_1(\om,\om_1)\,,\,K_2(\om,\om_1)$ and 
$\hat{\d}(\om)$ are multipliers on 
$C^*_r(G_A\times G_A)$. Using lemma \ref{K-del} 
we have the following sequence of 
equalities for any $\om,\om_1,\om_2,\om_3\in\sA(G_A)$
$$K_1(\om,\om_1)(\om_2\mt\om_3)=K_1(\om,\om_1)(\om_2\mt I)(I\mt\om_3)=
\hat{\d}(\om)(\om_2\mt\om_1)(I\mt\om_3)=$$
$$=\hat{\d}(\om)(I\mt\om_1)(\om_2\mt I)(I\mt\om_3)=
\hat{\d}(\om)(I\mt\om_1)(\om_2\mt \om_3).$$
Therefore,  multipliers $K_1(\om,\om_1)$ and $\hat{\d}(\om)(I\mt\om_1)$ agree on 
a dense subspace $\sA(G_A)\mt\sA(G_A)\subset C^*_r(G_A)\mt C^*_r(G_A)$, 
so they must be equal. The second equality
can be proved in the same way.\dowl\vs

\noindent
In this way the first part of the statement c) of thm. \ref{main} is proved.\vs

Now we will prove the density conditions. Comparing the formula for $K_1(f,f_1)$
with the expression for $W$ given in (\ref{W}) we see that:
$$K_1(f,f_1)(s,t)=k_1(s,t)(f\mt f_1)(W^{-1}(s,t))\,,\,
k_1(s,t):=|\Adab(\al(\br(s)^{-1}t))|^{-1}
|\Adbb(\br(s))|^{-1/2}$$
Because $k_1(s,t)\neq 0$,  we can apply lemma \ref{diff-W} to conclude that
elements in $\sA(G_A\times G_A)$ can be approximated  by elements of the 
form $K_1(\om,\om_1)$ in the norm given by $\om_0\mt\om_0$.  Then  from 
prop. \ref{r-m-prop} d) they can be approximated also in 
$C^*_r(G_A\times G_A)=C^*_r(G_A)\mt C^*_r(G_A)$ and since 
$\sA(G_A\times G_A)$ is dense $C^*_r(G_A\times G_A)$ we arrive at the 
conclusion: $span\{\Delta(a)(I\mt b)\,,\,a,b\in C^*_r(G_A)\}$ is dense in 
$C^*_r(G_A)\mt C^*_r(G_A)$.

The same argument can be applied to $K_2(f,f_1)$ since one easily checks that  
the mapping: 
$$G\times G\ni (s,t) \mapsto (\al(s)t,\bl(\al(s)t)^{-1}s)\in G\times G$$
is a  diffeomorphism of $G\times G$ (its inverse is given by: 
$(s,t)\mapsto(\bl(s)t,\al(\bl(s)t)^{-1})s)$).

This completes the proof of the theorem.\\
\dow

\noindent
In this way we have the main component of a quantum group structure on
$C^*_r(G_A)$ and in the next section we will look for the other ones.


\section{Hopf-like structure of $\sA(G_A)$.}
In this section we identify other mappings which appear
in the theory of quantum groups and show that the structure of 
algebra  $\sA(G_A)$ is very similar to the structure of a  Hopf algebra. 
Of course, because $\hat{\d}$ is not a mapping from 
$\sA(G_A)$ to $\sA(G_A)\mt\sA(G_A)$ and 
even not to multipliers of $\sA(G_A)\mt\sA(G_A)$ 
(in the sense of van Daele \cite{Vd}), the pair  $(\sA(G_A),\hat{\d})$ 
is not a true Hopf algebra.\vs

\noindent {\bf ``Coinverse''}\\
Let $Q$ be the function described in section 3, eq. (\ref{defQ}).
For $z\in \C$ we define a mapping $\tau_z$ by\notka{deftau}
\begin{equation}\label{deftau}
\tau_z:\sA(G_A)\ni \om \mapsto |Q|^{iz}\om\in \sA(G_A).
\end{equation}
Let us also define a mapping $R:\sA(G_A)\lra \sA(G_A)$ by\notka{defR}
\begin{equation}\label{defR}
R(\om)(g)(v\mt w):=\om(g^{-1})(w^{-1}\mt v^{-1})\,,\,\,
 {\rm where }\,\, v\in \lma T^l_gG_A\,,\, w\in \lma T^r_gG_A
\end{equation} 
and $v^{-1},w^{-1}$ denote the images of $v,w$ by the 
mapping $g\mapsto g^{-1}$. \\
Defining $Rf$ by $(Rf)\om_0:=R(f\om_0)$ 
we easily get  $(Rf)(g)=|\Adab (\al(g)\ar(g))|^{1/2}f(g^{-1})$.\\
Finally, let us  define $\kappa:\sA(G_A)\lra \sA(G_A)$ by \notka{def-kappa}  
\begin{equation}\label{def-kappa}
\kappa:=R \cdot \tau_{i/2}\,\, {\rm then}\,\, 
(\kappa f)(g)=|Q(g)|^{-1/2}|\Adab (\al(g)\ar(g))|^{1/2}f(g^{-1}).
\end{equation}
It is clear that all these mappings are linear.
\begin{prop}\label{kap-tau-R}\notka{kap-tau-R}
\begin{enumerate}
\item $\tau_z$ is a one-parameter (complex) group of automorphisms of 
 $\sA(G_A).$
\item $\tau_z \cdot *=*\cdot \tau_{\ove{z}}$
\item $R$ is an involutive $*$-antiautomorphism of $\sA(G_A)$.
\item $R\cdot \tau_z=\tau_z \cdot R$
\item $\kappa(\om_1\om_2)=\kappa(\om_2)\kappa(\om_1)$ 
\item $\kappa\cdot *\cdot\kappa\cdot *=id$.
\end{enumerate}
\end{prop}
{\em Proof:} Points 1. and 2. are 
immediate consequences of lemma \ref{koc}. They are true for any cocycle $Q$ 
on a differential groupoid \cite{PS2}.\vs\\ 
3. Since $R$ is implemented by a group inverse, it is clear that 
$R$ is involutive. From the formula for $Rf$ it immediately follows that 
$R$ commutes with *-operation on $\sA(G_A)$. 
To prove that $R$ is an antiautomorphism we compute:
$$R(f_1*f_2)(g)=|\Adab (\al(g)\ar(g))|^{1/2}(f_1*f_2)(g^{-1})=
|\Adab (a\ti{a})|^{1/2}\int_B\mu_L^2(b')f_1(\ti{a}^{-1}b') 
f_2(\bl(\ti{a}^{-1}b')^{-1}g^{-1}),$$
where $g=:a b=\ti{b}\ti{a}$.
On the other hand:
$$((Rf_2)*(Rf_1))(g)=\int_B \mu_L^2(b')(Rf_2)(ab') (Rf_1)(\bl(ab')^{-1}g)=$$
$$=\int_B \mu_L^2(b')f_1(g^{-1}\bl(ab')) f_2(b'^{-1}a^{-1})
|\Adab(a\ar(ab'))\Adab(\ar(ab')\ti{a})|^{1/2}=$$
$$=|\Adab(a\ti{a})|^{1/2}\int_B\mu_L^2(b')|\Adab(\ar(ab'))| 
f_1(g^{-1}\bl(ab')) f_2(b'^{-1}a^{-1})$$
Now apply the diffeomorphism $B\ni b'\mapsto \bl(\ti{a}b^{-1}b')\in B$ to this 
integral to get the previous one.\vs\\
4. This follows directly from the  definitions of $R$ and $\tau$.\vs\\
5 and 6. Are direct consequences of the definitions and the 
previous statements.
\\\dowl\vs

\noindent
{\bf ``Counit''.}\\
Since $F_l(e)=F_r(e)=B$ ($e$ is the neutral element  in $G$), 
for any $b\in B$ there  is  a mapping  
$\oml(b)\mt\omr(b)\ni\lambda\mt \rho \mapsto \lambda \rho \in \Om^1(T_bB)$.
In this way restriction of $\om\in\sA(G_A)$ to $B$ defines a smooth 
one-density on $B$ (with compact support). Integration of this density defines
a linear functional ${\epsilon}$ on $\sA(G_A): 
\sA(G_A)\ni \om\mapsto \int_B\om|_B\in C$. 
It is easy to see that
\begin{equation}\label{counit}
{\epsilon}(f\om_0)=\int_B\mu_L^2(b)|\Adbb(b)|^{1/2} f(b) .
\end{equation}
In the next proposition by $\kappa\mt id\,,\,id\mt\kappa\,,\,
\epsilon\mt id\,,\,id\mt \epsilon$ we 
denote the natural extensions of these mappings to 
$\sA(G_A\times G_A)$ and by $m:\sA(G_A\times G_A)\lra\sA(G_A)$ the natural 
extension of the multiplication map of $\sA(G_A)$.
\begin{prop}\label{hopf}\notka{hopf}
\begin{enumerate}
\item $\epsilon$ is a character on $\sA(G_A)$ 
\item $(\epsilon\mt id)\hat{\d}(\om)(I\mt\om_1)=\om \om_1$ and 
$(id \mt \epsilon)(\om_1\mt I)\hat{\d}(\om)=\om_1\om$;
\item $m(\kappa\mt id)\hat{\d}(\om)(I\mt\om_1)=\epsilon(\om)\om_1$ and 
$m(id\mt\kappa)(\om\mt I)\hat{\d}(\om_1)=\epsilon(\om_1)\om$.
\end{enumerate}
\end{prop}
{\em Proof:} 1. The definition of $\epsilon$ is a special case of the 
following one. Let $\G$ be a differential groupoid  and let $a\in E$ be such 
that $F_l(a)=F_r(a)$. Then we can define a linear functional $\epsilon_a$ 
on $\sA(\G)$ by $\epsilon_a(\om):=\int_{F_l(a)}\om|_{F_l(a)}$. We will prove
that  $\epsilon_a$ is a character on $\sA(\G)$.
 
Let $\om_1=f_1\om_0\,,\, \om_2=f_2\om_0$ then
$$\epsilon_a(\om_1\om_2)=
\int_{F_l(a)}\lo(x)\ro(x)\int_{F_l(a)}\lo^2(y)f_1(y)f_2(s(y)x)=
\int_{F_l(a)\times F_l(a)}\mu(x,y)f_1(y)f_2(s(y)x),$$ 
where  $\mu$ is a one density on $F_l(a)\times F_l(a)$ defined by
$$\mu(x,y)(v\mt w):=\lo(x)(v)\ro(x)(v)\lo^2(y)(w)\,,\,\,
v\in\lma T^l_x\G\,,\,w\in \lma T^l_y\G.$$
Now consider a diffeomorphism $\Phi$ of  $F_l(a)\times F_l(a)$ given by
$(x,y)\mapsto (s(y)x,y)$. Short computations show that
$(\Phi \mu)(z,t)=\nu(z,t)$ for 
$\nu(z,t)(v\mt w):=\lo(x)(v)\ro(x)(v)\lo(y)(w)\ro(y)(w)$. In this way our
integral is equal to
$$\int_{F_l(a)\times F_l(a)}\nu(x,y)f_1(x)f_2(y)=
\left(\int_{F_l(a)}\lo(x)\ro(x)f_1(x)\right)
\left(\int_{F_l(a)}\lo(y)\ro(y)f_2(y)\right)=
\epsilon_a(\om_1)\epsilon_a(\om_2).$$
To prove the second property, let us notice that the density 
$\ro(x)\lo(x)$ is real and invariant with respect to 
$s:F_l(a)\lra F_l(a)$. So we get
$$\epsilon_a(\om^*)=\int_{F_l(a)}f^*(x)\lo(x)\ro(x)=
\int_{F_l(a)}\ove{f(s(x))}\lo(x)\ro(x)=
\int_{F_l(a)}\ove{f(x)}\lo(x)\ro(x)=\ove{\epsilon_a(\om)}.$$ 
2. For $\om=F(\om_0\mt\om_0)$ we have 
$(\epsilon\mt id)(F(\om_0\mt\om_0))(g)=
(\int_B\mu_L^2(b) |\Adbb(b)|^{1/2} F(b,g) )\om_0(g).$
Using the prop. \ref{out-K} we compute:
$$(\epsilon\mt id)\hat{\d}(f \om_0)(I\mt f_1\om_0)(g)=
\left(\int_B\mu_L^2(b) |\Adbb(b)|^{1/2} (f*_\d(I\mt f_1))(b,g) \right)\om_0(g)=$$
$$=\left(\int_B\mu_L^2(b) 
|\Adab(\al(b^{-1}g))|^{-1}
f(b\al(b^{-1} g)) f_1(b^{-1}g)\right)\om_0(g)$$
Now apply a diffeomorphism $B\ni b\mapsto \bl(\al(g)^{-1}b)\in B$ to get 
the formula (\ref{product-Ga}) for $\om \om_1$.\vs\\
Since $(\om_1\mt I)\hat{\d}(\om)=(\hat{\d}(\om^*)(\om_1^*\mt I))^*$ we easily
get from prop. \ref{out-K}
$$(f_1\om_0\mt I)\hat{\d}(f\om_0)(g_1,g_2)=|\Adbb(\bl(g_2))|^{1/2}
f(\ar(\ar(g_1)\bl(g_2)^{-1})g_2)f_1(g_1\bl(g_2)^{-1})(\om_0\mt\om_0)(g_1,g_2)$$
and  we compute
$$(id\mt\epsilon)(f_1\om_0\mt I)\hat{\d}(f \om_0)(g)=
\left(\int_B\mu_L^2(b) |\Adbb(b)|^{1/2} ((f_1\mt I)*_\d f)(g,b)\right)
\om_0(g)=$$
$$=\left(\int_B\mu_L^2(b) |\Adbb(b)| f_1(gb^{-1}) f(\ar(\ar(g)b^{-1})b)\right)
\om_0(g)$$
Again, applying a diffeomorphism $B\ni b\mapsto \br(g)b^{-1}\in B$ we arrive 
at the formula (\ref{product-Ga}) for $\om_1 \om$.\vs\\
3.  For $F(\om_0\mt\om_0)\in\sA(G_A\times G_A)$ the mapping $m$ is given by
(compare (\ref{product-Ga}))
$$m(F(\om_0\mt\om_0))=(mF)\om_0\,\,{\rm  and}\,\, 
(mF)(g):=\int_B\mu_L^2(b) F(\al(g)b\,,\bl(\al(g)b)^{-1}g)$$
From the formula for $\kappa$ (\ref{def-kappa}) and prop. \ref{out-K} 
we get
$$(\kappa\mt id)\hat{\d}(f\om_0)(I\mt f_1\om_0)(g_1,g_2)=
(\kappa\mt id)(f*_\d(I\mt f_1))(g_1,g_2)
(\om_0\mt\om_0)(g_1,g_2)\,\,{\rm and}$$
$$(\kappa\mt id)(f*_\d(I\mt f_1))(g_1,g_2)=
|Q(g_1)|^{-1/2}|\Adab (\al(g_1)\ar(g_1))|^{1/2} 
(f*_\d(I\mt f_1))(g_1^{-1},g_2).$$
Let $g=a b$, we need a   value of this function for $g_1=a\ti{b}$ and 
$g_2=\bl(a\ti{b})^{-1}a b$. For this points  we have
$\br(g_1^{-1})=\bl(a\ti{b})^{-1}$ and 
$$(\kappa\mt id)(f*_\d(I\mt f_1))(a\ti{b},\bl(a\ti{b})^{-1}a b)=$$
$$=|Q(a\ti{b})|^{-1/2}|\Adab (a \ar(a\ti{b}))|^{1/2}
(f*_\d(I\mt f_1))(\ti{b}^{-1}a^{-1},\bl(a\ti{b})^{-1}a b)=$$
$$=|Q(a\ti{b})|^{-1/2}|\Adab (a \ar(a\ti{b}))|^{1/2}
|\Adab(a))|^{-1} |\Adbb(\bl(a\ti{b}))|^{1/2}
f(\ti{b}^{-1}) f_1(a b)=$$
$$=|\Adbb(\ti{b})|^{1/2}
f(\ti{b}^{-1})  f_1(a b),$$
where  the last equality follows from (\ref{Q-mod}).
In this way
$$m(\kappa\mt id)(f*_\d(I\mt f_1))(g)=f_1(g) \int_B\mu_L^2(b)|\Adbb(b)|^{1/2}
f(b^{-1})=\epsilon(f\om_0)f_1(g).$$
The second equality can be proved by similar computations.\\
\dowl
\begin{re}{\em 
The first two statements of the above proposition can also be proved by
the following observation. If in a differential groupoid $\G$ there is
$a\in E$ such that $F_l(a)=F_r(a)$, then a relation $h_a:\G\rel \{1\}$ 
defined by $Gr(h_a):=\{(1,g) : g\in F_l(a)\}$ is a morphism of differential 
groupoids. Since $\sA(\{1\})=C$ the associated mapping 
$\widehat{h_a}$ takes values in $C$
and the functional $\epsilon_a$ is just $\widehat{h_a}$. For the functional
$\epsilon$,  $a=e$ (the neutral element) and the relation is $e_B^T$ 
($e_B$-identity of $G_B$). So the second statement is implied by 
$(e_B^T\times id)\d=(id\times e_B^T)\d=id$, and this is just the transposition
of the equality: $m_B(e_B\times id)=m_B(id\times e_B)=id$.}
\end{re}
The statements in the last proposition are just axioms for a Hopf algebra
rewritten in a way that makes sense in our situation. Whether this need be
formalized is a question of finding other interesting examples.

Now we describe, how $\hat{\delta}$ commutes with $R$ and $\tau_z$.
Let us define a  mapping $\sim: \sA(G_A\times G_A)\lra\sA(G_A\times G_A)$ by:
$(\sim\om)(g_1,g_2):=\om(g_2,g_1)$.
\begin{prop}\label{delta-R-tau}\notka{delta-R-tau} 
For any $z\in\C$ and $\om,\om_1,\om_2\in\sA(G_A)$: \notka{delta-R}
\begin{eqnarray}\label{delta-R}
\hat{\d}(R\om)(R\om_1\mt R\om_2)& = &
\sim (R\mt R)((\om_2\mt\om_1)\hat{\d}(\om))\\
\label{delta-tau}
\hat{\d}(\tau_z\om)(\om_1\mt\om_2)& =&
(\tau_z\mt\tau_z)(\hat{\d}(\om)(\tau_{-z}\om_1\mt\tau_{-z}\om_2))
\end{eqnarray}\notka{delta-tau}\vsm\vsm\vsm\vsm\vsm
\end{prop}
{\em Proof:} Let $\om=:f\om_0\,,\,\om_1=:f_1\om_0\,,\,\om_2=:f_2\om_0$. 
Using the same notation as before, formulae for $\hat{\d}$ (\ref{komno})
and $R$, let us compute the left hand side of the first equality
$$(Rf)*_\d(Rf_1\mt Rf_2)(a_1 b_1, a_2 b_2)=$$
$$=\int_{B} \mu_L^2(b) |\Adbb(\bl(a_2 b))|^{-1/2}
(Rf)(a_1a_2b) (Rf_1)(\bl^{-1}(a_1a_2 b) a_1 b_1)(Rf_2)
(\ar(a_2 b) b^{-1} b_2)=$$
$$=\int_{B} \mu_L^2(b) |\Adbb(\bl(a_2 b))|^{-1/2}
|\Adab (a_1 a_2\ar(a_1a_2b))|^{1/2}f(b^{-1}a_2^{-1}a_1^{-1})\times$$
$$\times |\Adab (\ar(a_1a_2b)\al(b^{-1}a_2^{-1})\ar(a_1b_1))|^{1/2}
f_1(b_1^{-1} a_1^{-1}\bl(a_1a_2 b))\times$$
$$\times |\Adab (\ar(a_2b)\ar(a_2 b_2))|^{1/2}
f_2(b_2^{-1}b\ar(a_2 b)^{-1})=$$
$$=\int_{B} \mu_L^2(b) |\Adbb(\bl(a_2 b))|^{-1/2}
|\Adab (a_1 a_2\ar(a_1a_2b)^2 \ar(a_1b_1) \ar(a_2 b_2))|^{1/2}\times$$
$$\times f(b^{-1}a_2^{-1}a_1^{-1})
f_1(b_1^{-1} a_1^{-1}\bl(a_1a_2 b))
f_2(b_2^{-1}b\ar(a_2 b)^{-1})=$$
$$=|\Adab (a_1 a_2\ar(a_1b_1) \ar(a_2 b_2))|^{1/2}\times$$
$$\times \int_{B} \mu_L^2(b) |\Adbb(\bl(a_2 b))|^{-1/2}
|\Adab (\ar(a_1a_2b))| f(b^{-1}a_2^{-1}a_1^{-1})
f_1(b_1^{-1} a_1^{-1}\bl(a_1a_2 b))
f_2(b_2^{-1}b\ar(a_2 b)^{-1})$$
And the right hand side
$$\sim (R\mt R)((f_2\mt f_1)*_\d f)(a_1 b_1,a_2 b_2)=
(R\mt R)((f_2\mt f_1)*_\d f)(a_2 b_2,a_1 b_1)=$$
$$=|\Adab (a_2 a_1 \ar(a_2 b_2)\ar(a_1 b_1))|^{1/2}
(f_2\mt f_1)*_\d f)(b_2^{-1} a_2^{-1},b_1^{-1} a_1^{-1})=$$
$$=|\Adab (a_2 a_1 \ar(a_2 b_2)\ar(a_1 b_1))|^{1/2}
\ove{(f^**_\d(f_2^*\mt f_1^*))(s_B(a_2 b_2),s_B(a_1 b_1))}=$$
$$=|\Adab (a_2 a_1 \ar(a_2 b_2)\ar(a_1 b_1))|^{1/2}
\ove{(f^**_\d(f_2^*\mt f_1^*))(a_2^{-1}\bl(a_2 b_2),a_1^{-1}\bl(a_1 b_1))}=$$
$$=|\Adab (a_1 a_2 \ar(a_1 b_1)\ar(a_2 b_2))|^{1/2}\times$$
$$\times \ove{\int_{B} \mu_L^2(b) |\Adbb(\bl(a_1^{-1} b))|^{-1/2}
f^*(a_2^{-1} a_1^{-1} b) 
f_2^*(\bl^{-1}(a_2^{-1}a_1^{-1} b)a_2^{-1}\bl(a_2 b_2))
f_1^*(\ar(a_1^{-1} b) b^{-1}\bl(a_1 b_1) )}=$$
$$=|\Adab (a_1 a_2 \ar(a_1 b_1)\ar(a_2 b_2))|^{1/2}\times$$
$$\times \int_B \mu_L^2(b) |\Adbb(\bl(a_1^{-1} b))|^{-1/2}
f(\bl(a_2^{-1} a_1^{-1} b)^{-1}a_2^{-1}a_1^{-1}) f_1(b_1^{-1}a_1^{-1}b)
f_2(b_2^{-1}a_2^{-1}\br(b^{-1}a_1)^{-1})$$
Now use a diffeomorphism $b\mapsto \bl(a_2^{-1}a_1^{-1}b)$ to get the equality
of two integrals.\vs

To prove the second equality one immediately verifies that it is 
enough to prove that:
$$Q(a_1a_2b)=\frac{Q(a_1b_1)Q(a_2b_2)}
{Q(\bl(a_1a_2b)^{-1}a_1b_1)Q(\ar(a_2 b)b^{-1}b_2)}$$
But since $Q$ is a cocycle on $G_A$ and $G_B$ we have equalities:
$Q(a_1a_2b)=Q(a_1\bl(a_2b))Q(a_2b)$,\\ 
$Q(\bl(a_1a_2b)^{-1}a_1b_1)=Q(\bl(a_1a_2b)^{-1}a_1)Q(a_1b_1)$ 
 and  $Q(a_1\bl(a_2b))=Q(\bl(a_1a_2b)^{-1}a_1)^{-1}$. Using them
we easily arrive at  the desired result.\\\dowl\vsm\vsm\\


\section{Haar measure.}

We begin by recalling some results from \cite{PS2}.
Let $(\G,m,s,E)$ be a differential groupoid. Choose some real, non-vanishing,
right-invariant half-density $\ti{\ro}$. Let $\ti{\om}_0=\ti{\lo}\mt\ti{\ro}$ 
be the corresponding bidensity. Let us also choose some real, non-vanishing 
half-density $\nu$ on $E$.
For such a pair we define:
\begin{itemize}
\item A smooth function $\sigma: \G\lra R\setminus\{0\}$ defined by  
$\sigma(g)(\ti{\ro}\mt\nu)(g)=(\ti{\lo}\mt\nu)(g)$. We call $\sigma$ 
{\em the modular function} associated with a pair $(\ti{\ro},\nu)$.
This function is a  one cocycle on $\G$ (i.e. for any composable $g_1,g_2\in\G$: 
$\sigma(g_1g_2)=\sigma(g_1)\sigma(g_2)$). Therefore, it defines a one 
parameter (complex) group of automorphisms (not $*$-automorphisms!) of 
$\sA(\G): \sigma_z(\om)(g):=|\sigma(g)|^{2 i z}\om(g)$;  
\item A positive linear functional $h$ on $\sA(\G)$: $h(f\ti{\om_0}):=\int_E \nu^2 f$.
\item A linear mapping 
$\hat{h}:\sA(\G)\ni f \ti{\om}_0\mapsto  f \ti{\ro}\mt\nu \in L^2(\G)$.
\end{itemize}
The next proposition describes basic properties of these objects.
\begin{prop}\label{GNS}\notka{GNS}
{\rm \cite{PS2}} For any $\om,\om_1\in\sA(\G)$:
\begin{itemize} 
\item  $ h(\om^*\om)=(\hat{h}(\om)\,|\,\hat{h}(\om))$ 
(the right hand side is a scalar product in $L^2(\G)$);  
\item $h(\sigma_z(\om)\om_1)=h(\om_1\sigma_{z-i}(\om))$;
\item $\hat{h}(\om\om_1)=\pi_{id}(\om)\hat{h}(\om_1)$
\end{itemize}\dowl
\end{prop}

Now we specify the choice of $\ti{\ro}$ and $\nu$ for our groupoid $G_A$.
Let $\mu_0\neq 0$ be  a real, half-density on $T_eB$ and let us define 
a right-invariant half-density $\ti{\ro}$ on $G_A$ by the formula
$$\ti{\ro}(g)(w):=\mu_0(wg^{-1})\,,\,w\in \lma T^r_gG_A.$$
Let $\ti{\lo}$ be the  corresponding left-invariant half-density  and 
$\ti{\om}_0=\ti{\lo}\mt\ti{\ro}$. Short computations show that
$\ti{\om}_0(g)=|\Adab(\al(g)\ar(g))|^{1/2}\om_0$ 
(where $\om_0$ denotes our standard bidensity
constructed from $\mu_0$ as in (\ref{lambda},\ref{ro})).
Let us also choose some $\nu_0\neq 0$ -- a real, half-density on $T_eA$ and 
let $\nu_r$ be the corresponding right-invariant half-density on $A$.
In  our standard representation of bidensities we get \notka{haar}
$$\label{haar}h(f\om_0)=\int_A \nu^2_r(a)|\Adab(a)|^{-1} f(a).$$
From the definitions of $h$, $\sigma$ and $\tau$ it is clear that
$h\tau_z=h$ and $\sigma_z\tau_w=\tau_w\sigma_z$.
In the next lemma we give the formula for  the associated modular function.
\begin{lem}\label{modular}\notka{modular} Let $\ti{\ro}$ and $\nu_r$ 
be as above. Then the modular function $\sigma$ is given by
$$\sigma(g)=\left|\frac{\Adaa(\al(g))\Adba(\br(g))}
{\Adbb(\bl(g))\Adaa(\ar(g))}\right|^{1/2}.$$
\end{lem}
{\em Proof:} From the definition of $\ti{\lo}$ one can easily see that 
$\ti{\lo}(g)(g v)=|\Adab(\ar(g))|^{1/2} \mu_0(v)\,,\,v\in \lma T_eB$.
The  subspace $gT_eA\subset T_gG$ is complementary to 
$T^r_gG_A$ and $T^l_gG_A$ i.e. $T_gG=gT_eA\ms T^r_gG_A=gT_eA\ms T^l_gG_A$. 
Therefore,  for $u\in \lma T_eA\,,\, v\in\lma T_eB$ we have the equalities:
$$(\ti{\ro}\mt\nu_r)(g)(vg\wedge gu)=\mu_0(v)\nu_r(\ar(g))(\ar(gu))=
|\Adaa(\ar(g))|^{1/2}\mu_0(v)\nu_0(u)$$
$$(\ti{\lo}\mt\nu_r))(g)(g v\wedge g u)=|\Adab(\ar(g))|^{1/2} \mu_0(v)
\nu_r(\al(g))(\al(g u))=$$
$$=|\Adab(\ar(g))|^{1/2}|\Adba(\br(g))|^{1/2}|\Adaa(\al(g))|^{1/2}
\mu_0(v)\nu_0(u)$$
Now, since $vg\wedge gu=g(g^{-1}vg)\wedge gu=
\det(P_B Ad(g^{-1}) P_B) (g v\wedge g u)$ we get the expression for $\sigma(g)$
$$\sigma(g)=|\Adab(\ar(g))|^{1/2}|\Adba(\br(g))|^{1/2}|\Adaa(\al(g))|^{1/2}
|\det(P_B Ad(g^{-1}) P_B)|^{1/2}|\Adaa(\ar(g))|^{-1/2},$$
and the result follows from the equality 
$\det(P_B Ad(g^{-1}) P_B)=\Adab(\ar(g))^{-1}\Adbb(\bl(g))^{-1}$. 
\\\dowl\vs

Let $h\mt id: \sA(G_A\times G_A)\lra \sA(G_A)$ be the natural extension
of $h\mt id: \sA(G_A)\mt\sA(G_A)\lra \sA(G_A)$ i.e.
$((h\mt id)(F (\om_0\mt\om_0)))(g):=
\left(\int_A |\Adab(a)|^{-1}\nu^2_r(a) F(a,g)\right)\om_0(g)$.
The next proposition states that the functional $h$ is ``right-invariant'', 
relates $h$ to $\kappa$ and describes commutation of $\hat{\d}$ with the 
modular group $\sigma_z$. Let us recall that a functional $\varphi$ on
a Hopf algebra $(\sA,\Delta)$ is right-invariant iff 
$(\varphi\mt id)\Delta(a)=\varphi(a) I$.

\begin{prop}\label{proteza}\notka{proteza} 
For any  $\om\,,\,\om_1\,,\,\om_2\in\sA(G_A)$:
\begin{eqnarray}
(h\mt id)(\hat{\d}(\om_1)(I\mt \om_2))& = & h(\om_1)\om_2,\\
\label{kappa-h}(h\mt id)(\hat{\delta}(\om_1^*)(\om_2\mt I)) & = &
\kappa((h\mt id)((\om_1^*\mt I)\hat{\delta}(\om_2))\\
\hat{\d}(\sigma_z(\om))(\om_1\mt\om_2)&=&
(\sigma_z\mt\tau_z)(\hat{\d}(\om)(\sigma_{-z}\om_1\mt\tau_{-z}\om_2)
\end{eqnarray}
\end{prop}\vsm\vsm
{\em Proof:} Using prop. \ref{out-K} and right invariance of $\nu_r$ we get
$$(h\mt id)(\hat{\d}(f_1\om_0)(I\mt f_2\om_0))(g)=
\left(\int_A\nu^2_r(a)|\Adab(a)|^{-1}|\Adab(\al(g))|^{-1}
f_1(a\al(g))f_2(g)\right)\om_0(g)=$$
$$=\left(\int_A\nu^2_r(a)|\Adab(a\al(g))|^{-1}
f_1(a\al(g))\right)\om_2(g)=h(\om_1)\om_2(g).$$
To prove the second formula, we again use prop. \ref{out-K}, right invariance
of $\nu_r$ and the formula for 
$(\om_1^*\mt I)\hat{\delta}(\om_2)$ given in the 
proof of the proposition \ref{hopf} to obtain the left hand side of the equality 
(\ref{kappa-h}).
$$(h\mt id)(\hat{\delta}((f_1\om_0)^*)(f_2\om_0\mt I))(g)=
\left(\int_A\nu_r^2(a) |\Adab(a)|^{-1}(f_1^**_\delta(f_2\mt I))(a,g)\right)
\om_0(g)=$$
$$=\left(|\Adbb(\bl(g))|^{-1/2}\int_A\nu_r^2(a)|\Adab(a)|^{-1}
\overline{f_1(s_A(ag))}f_2(\bl(ag)^{-1}a)\right)\om_0(g)$$
And the right hand side of this equality:
$$\kappa((h\mt id)(((f_1\om_0)^*\mt I)\hat{\delta}(f_2\om_0))(g)=$$
$$=\left(|Q(g)|^{-1/2}|\Adab(\al(g)\ar(g))|^{1/2}\int_A\nu_r^2(a)
|\Adab(a)|^{-1}((f_1^*\mt I)*_\delta f_2)(a,g^{-1})\right)\om_0(g)=$$
$$=\left(|Q(g)|^{-1/2}|\Adab(\al(g)\ar(g))|^{1/2}|\Adbb(g^{-1})|^{1/2}
\right.\times $$
$$\left.\times\int_A\nu_r^2(a)|\Adab(a)|^{-1}
\overline{f_1(s_A(a\br(g)))}f_2(\ar(a\br(g))g^{-1})\right)\om_0(g) $$
Now we use the relationship among $Q$ and modular functions (\ref{Q-mod}) to
convert this expression into
$$\left(|\Adab(\al(g))|^{1/2}|\Adbb(\bl(g))|^{-1/2}
\int_A\nu_r^2(a)|\Adab(a)|^{-1}\overline{f_1(s_A(a\br(g)))}
f_2(\ar(a\br(g))g^{-1})\right)\om_0(g)$$
Since $a\br(g)=a\al(g)^{-1}g$, 
$\ar(a\br(g))g^{-1}=\bl(a\al(g)^{-1}g)^{-1}a\al(g)^{-1}$ 
and $\nu_r$ is right-invariant this 
is equal to the left hand side.\\
From  the formula for $\hat{\d}$ (\ref{komno}) and definitions of $\sigma$
and $\tau$ it easily follows that to prove the third statement
it is sufficient to show the equality
$$\frac{\Adaa(a_1a_2)\Adba(b)}{\Adaa(\ar(a_1a_2b))\Adbb(\bl(a_1a_2b))}=$$
$$=\frac{\Adaa(a_1)\Adba(b_1)Q(a_2b_2)\Adbb(\bl(\bl^{-1}(a_1a_2b)a_1b_1))
\Adaa(\ar(\bl^{-1}(a_1a_2b)a_1b_1))}
{\Adbb(\bl(a_1b_1))\Adaa(\ar(a_1b_1))\Adaa(\al(\bl^{-1}(a_1a_2b)a_1b_1))
\Adba(\br(\bl^{-1}(a_1a_2b)a_1b_1)) Q(\ar(a_2b)b^{-1}b_2)}$$
and this is straightforward because of  equalities:
$$\bl(\bl^{-1}(a_1a_2b)a_1b_1)=\bl^{-1}(a_1a_2b)\bl(a_1b_1)\,,\,
\ar(\bl^{-1}(a_1a_2b)a_1b_1)=\ar(a_1b_1),$$
$$\al(\bl^{-1}(a_1a_2b)a_1b_1)=\ar(a_1a_2b)\ar^{-1}(a_2b)\,,\,
\br(\bl^{-1}(a_1a_2b)a_1b_1)=\bl^{-1}(a_2b)b_1,$$
$$Q(a_2b_2)=Q(\ar(a_2b)b^{-1}b_2)Q(a_2b),$$ and equation (\ref{Q-mod}).\\
\dowl \\
The relationships among $h$, $\hat{h}$ and $\sigma_t$ given in prop \ref{GNS}
suggest that $\hat{h}$ can be extended to a GNS mapping and $h$ to a KMS weight
on $C^*_r(G_A)$. And indeed this is true and was proved in \cite{PS2} for
a general differential groupoid.
\begin{prop}{\rm \cite{PS2}} The mapping $\hat{h}$ is closable and 
defines GNS mapping from $C^*_r(G_A)$ to $L^2(G)$. Consequently, the linear
functional $h$ can be extended to a densely defined, lower semi-continuous
weight on $C^*_r(G_A)$, which is a KMS weight with the modular group $\sigma_t$.
\end{prop}

In the remaining part of this section we show that this weight is right-invariant. 
The crucial (although not difficult) step is lemma 
\ref{W-h}, and then we use results from \cite{Vaes-Kust} and \cite{WMN}.

Let  $H$ be a Hilbert space. For two vectors $x,y\in H$, we define a linear 
functional  $\eta_{xy}$ on $B(H)$ by $\eta_{xy}(a):=(x|a y)$. Then for
$a\in B(H\mt H)$ we define operators $(id\mt\eta_{xy})a$ and 
$(\eta_{xy}\mt id)a$ by  $(z|((id\mt\eta_{xy})a)t):=(z\mt x|a(t\mt y))$ and
$(z|((\eta_{xy}\mt id)a)t):=(x\mt z|a(y\mt t))$.
\begin{lem}\label{W-h}\notka{W-h}
Let $\hat{h}:C^*_r(G_A)\lra L^2(G)$ be the GNS-mapping associated with the weight
$h$ and $W$ the multiplicative unitary defined in (\ref{W}). 
For $x,y\in L^2(G)$ and  $a\in D(\hat{h})$ an element 
$\hat{h}((id\mt\eta_{xy})\pi_{\delta}(a))$ belongs to $D(\hat{h})$ and 
$$\hat{h}((id\mt\eta_{xy})\pi_{\delta}(a))=((id\mt\eta_{xy})W)\hat{h}(a)$$
\end{lem}
{\em Proof: } First, let us notice that since $\hat{h}$ is a closed 
mapping and $\sA(G_A)$ is its  core, it is enough to prove the formula for
$x,y\in \sD(G)$ and $a\in \sA(G_A)$. Indeed, let  $a\in D(\hat{h})$ and 
$x=\hat{h}(a)$. This means that 
$a=\lim\om_n\,,\,\om_n\in \sA(G_A)$ and $x=\lim \hat{h}(\om_n)$.
It is shown below that  $(id\mt\eta_{xy})\pi_{\delta}(\om_n)$ is a sequence
in $\sA(G_A)$ and it  converges to  $(id\mt\eta_{xy})\pi_{\delta}(a)$. 
So  $\hat{h}((id\mt\eta_{xy})\pi_{\delta}(\om_n))=
((id\mt\eta_{xy})W)\hat{h}(\om_n)$. But since $\lim\hat{h}(\om_n)=x$, the
sequence on the right hand side is convergent to $((id\mt\eta_{xy})W)x$. 
Therefore $(id\mt\eta_{xy})\pi_{\delta}(a)$ is in $D(\hat{h})$ and 
$\hat{h}((id\mt\eta_{xy})\pi_{\delta}(a))=((id\mt\eta_{xy})W)\hat{h}(a)$.
Since $\sD(G)$ is dense in $L^2(G)$, similar arguments show that it is enough
to check the equality for $x,y\in \sD(G)$.

Now we compute $(id\mt\eta_{xy})\pi_{\delta}(\om)$ for 
$x,y\in \sD(G)$ and $\om\in\sA(G_A)$.
Let  $\psi_0:=\ro\mt\nu_r$, where $\ro$ is our standard right-invariant 
half-density defined in (\ref{ro}) and $\nu_r$ is as  in the definition of $h$. 
We  write: $x=:f_x\psi_0,\, y=:f_y\psi_0,\,  \om=:f\om_0$. With this notation we get
$(id\mt\eta_{xy})\pi_{\delta}(\om)=f_{xy}\om_0$, where
\notka{slice}
\begin{equation}
\label{slice}f_{xy}(ab):=|\Adbb(b))|^{-1/2}
\int_{G}\psi_0^2(g_2)\ove{f_x(g_2)}|\Adab(\ar(a_2^{-1} b))|
f(aa_2\bl(a_2^{-1}b))f_y(b^{-1}g_2).
\end{equation}

To prove  this formula it is enough to show that 
$(z|\pi_{id}(f_{xy}\om_0)t)=(z\mt x|\pi_{\delta}(\om)(t\mt y))$ 
for any $z,t\in \sD(G)$. 
Let $z=:f_z\psi_0, t=:f_t\psi_0.$
From the formula for 
$\hat{\delta}$ (\ref{komno}) we get
$$(z\mt x|\pi_{\delta}(\om)(t\mt y))=
\int_{G\times G}(\psi_0\mt\psi_0)^2(g_1,g_2)
\ove{f_z(g_1)}\ove{f_x(g_2)}\times$$
$$\times\int_B\mu_L^2(b)|\Adbb(\bl(a_2b))|^{-1/2} 
f(a_1a_2b)f_t(\bl(a_1a_2b)^{-1}a_1b_1)f_y(\ar(a_2b)b^{-1}b_2)=$$
$$=\int_{G}\psi_0^2(g_1)\ove{f_z(g_1)}\int_{G}\psi_0^2(g_2)\ove{f_x(g_2)}
\times$$
$$\times\int_B\mu_L^2(b)|\Adbb(\bl(a_2b))|^{-1/2} 
f(a_1a_2b)f_t(\bl(a_1a_2b)^{-1}a_1b_1)f_y(\ar(a_2b)b^{-1}b_2),$$
where $g_1=a_1b_1, \,g_2=a_2 b_2.$\\
Applying a diffeomorphism $B\ni b\mapsto \bl(a_2 b)\in B$, we can rewrite
integral over $B$ as 
$$\int_B\mu_L^2(b)|\Adbb(b))|^{-1/2}|\Adab(\ar(a_1^{-1} b))|
f(a_1a_2\bl(a_2^{-1}b))f_y(b^{-1}a_2 b_2)f_t(\bl(a_1b)^{-1}a_1b_1),$$
and interchanging the order of integration we get
$$(z\mt x|\pi_{\delta}(\om)(t\mt y))=
\int_{G}\psi_0^2(g_1)\ove{f_z(g_1)}\int_B\mu_L^2(b)|\Adbb(b))|^{-1/2}\times$$
$$\times\left[\int_{G}\psi_0^2(g_2)\ove{f_x(g_2)}|\Adab(\ar(a_1^{-1} b))|
f(a_1a_2\bl(a_2^{-1}b))f_y(b^{-1}g_2)\right]f_t(\bl(a_1b)^{-1}g_1).$$
Comparing this with the formula for a multiplication in $G_A$ (\ref{product-Ga}) 
and the definition of $f_{xy}$ we obtain the result.

To prove the lemma it remains to show that
$(z|\hat{h}(f_{xy}\om_0))=(z\mt x|W(\hat{h}(\om)\mt y))$.
Using  the definition of $\hat{h}$ we obtain $\hat{h}(f\om_0)=\hat{f}\psi_0$, for
$\hat{f}(g):=f(g)|\Adab(\ar(g))|^{-1/2}$.
Let us compute the right hand side of our equality
$$(z\mt x|W(\hat{h}(\om)\mt y))=
\int_{G\times G}(\psi_0\mt\psi_0)^2(g_1,g_2)\overline{f_z(g_1)}
\overline{f_x(g_2)}[W(\hat{f}\mt f_y)](g_1,g_2)=$$
$$=\int_G\psi_0^2(g_1)\overline{f_z(g_1)}\int_G\psi_0^2(g_2)
\overline{f_x(g_2)}[W(\hat{f}\mt f_y)](g_1,g_2)$$
In this way we have to prove the equality
$$f_{xy}(g)|\Adab(\ar(g))|^{-1/2}=\int_G\psi_0^2(g_2)
\overline{f_x(g_2)}[W(\hat{f}\mt f_y)](g,g_2),\qquad\qquad(**)$$
where  $W(f\psi_0\mt k\psi_0)=:W(f \mt k)(\psi_0\mt\psi_0)$.\\
The formula for the action of $W$ was given after the  lemma \ref{W'}. 
Notice, that now $\psi_0$ is not the same as $\psi_0$
defined before that lemma, so the expressions for 
$W(\hat{f}\mt f_y)$ are different. But since both $\psi_0$'s 
can be easily compared, we immediately deduce needed formula from the one
contained in lemma \ref{W'}:\vsm
$$[W(\hat{f}\mt f_y)](g,g_2)=|\Adbb(\br(g))\Adab(\al(\br(g)^{-1}g_2))|^{-1/2}
\hat{f}(g\al(\br(g)^{-1}g_2))f_y(\br(g)^{-1}g_2)$$
To verify the equality $(**)$ use this expression,  expression 
for $f_{xy}$ (\ref{slice}) and $\hat{f}$. This completes the proof of the 
lemma.\\
\dowl\\
Now we can prove right invariance of $h$ with respect to positive
vector functionals in a GNS space for $h$. This immediately follows from
\begin{prop}
For any $a\in D(\hat{h})$ and $x\in L^2(G)$:
$$h((id\mt\eta_{xx})\pi_\delta(a^*a))=||x||^2h(a^*a).$$
\end{prop}
{\em Proof: }
Let $e_i$ be an orthonormal basis in $L^2(G)$. For $x\in L^2(G)$ let
$$B_n:=\sum_{i=1}^n[(id\mt\eta_{e_ix})\pi_\delta(a)]^*
[(id\mt\eta_{e_ix})\pi_\delta(a)].$$
Then $B_n\in C^*_r(G_A)$ and $B_n$ converges strictly to 
$(id\mt \eta_{xx})\pi_\delta(a^*a)$ (see \cite{Vaes-Kust}).
On the other hand
$$h(B_n)=\sum_{i=1}^n 
\left(\hat{h}[(id\mt\eta_{e_ix})\pi_\delta(a)]|\hat{h}
[(id\mt\eta_{e_ix})\pi_\delta(a)]\right)=
\sum_{i=1}^n\left([(id\mt\eta_{e_ix})W]\hat{h}(a)|
[(id\mt\eta_{e_ix})W]\hat{h}(a)\right)=$$
$$=(\hat{h}(a)|\sum_{i=1}^n[(id\mt\eta_{e_ix})W]^*[(id\mt\eta_{e_ix})W]\hat{h}(a)).$$
Because $W\in M(CB(L^2(G))\mt C^*_r(G_A))$ 
(see remarks after lemma \ref{W-bis}),  the sequence \\
$\sum_{i=1}^n[(id\mt\eta_{e_ix})W]^*[(id\mt\eta_{e_ix})W]$ 
converges  strongly to $(id\mt\eta_{xx})(W^*W)=||x||^2 I$. Therefore 
$h(B_n)$ converges to $||x||^2 h(a^*a)$. Since the weight $h$ is strictly lower
semi-continuous we get the result.\\
\dow

\noindent
The last step is to prove that this is enough for full right invariance. 
To this end we use the following
\begin{tw}{\rm \cite{WMN}}
Let $h$ be a densely defined, lower semicontinuous weight on a separable 
$C^*$-algebra $A$ and  $(H,\eta_h,\pi_h)$ be a GNS-triple.
There exists a sequence of vectors $\Omega_n\in H$ such that:
\begin{enumerate}
\item $\pi_h(A)\Omega_n$ and $\pi_h(A)\Omega_m$ are orthogonal for $m\neq n$;
\item $a\in D(\eta_h) \iff \sum ||\pi_h(a)\Omega_n||^2 < \infty$;
\item $\eta_h(a)=\sum \pi_h(a)\Omega_n$;
\item $H=\oplus H_n\,,\,H_n:=\overline{\pi_h(A)\Omega_n}$
\end{enumerate}
\end{tw}\dow

To simplify notation we put:  $\eta_i:=\eta_{\Omega_i\Omega_i}$  and, for  
$a,b\in D(\hat{h})$,  $\eta_{i a}:=\eta_{\Om_i\hat{h}(a)}$,\, 
$\eta_{ab}:=\eta_{\hat{h}(a)\hat{h}(b)}$. We also denote $C^*_r(G_A)$ by $A$.
From lemma \ref{W-h} and the theorem above we easily obtain three useful
formulae:
$$(\eta_i\mt id)\Delta(a)=(\eta_{ia}\mt id)W,\,\,a\in D(\hat{h}),$$
$$\Delta(a)(\Om_n\mt x)=(p_n\mt I)W(\hat{h}(a)\mt x),\,\,a\in D(\hat{h}),$$
$$(\eta_{i(a^*a)}\mt id)W=(\eta_a\mt id)(W^*(p_i\mt I)W),$$
where $p_n$ denotes projection onto the closure of $A\Om_n$.\\
Indeed, let $x,y,z\in L^2(G)$ then
$$(x|(\eta_{ia}\mt id)Wy)=
(\Om_i\mt x|W(\hat{h}(a)\mt y))=(\Om_i|(id\mt\eta_{xy})W\hat{h}(a))=$$
$$=(\Om_i|\hat{h}((id\mt\eta_{xy})\Delta(a)))=
(\Om_i|(id\mt\eta_{xy})\Delta(a)\Om_i)=(x|(\eta_i\mt id)\Delta(a)y),$$
and this proves the first one. To prove the second one,  let us compute
$$(x\mt z| W(\hat{h}(a)\mt y))=(x|(id\mt \eta_{zy})W\hat{h}(a))=
(x|\hat{h}((id\mt \eta_{zy})\Delta(a)))=
\sum_n(x|(id\mt\eta_{zy})\Delta(a)\Om_n),$$
so
$$(x\mt z|(p_n\mt I) W(\hat{h}(a)\mt y))=(x|(id\mt\eta_{zy})\Delta(a)\Om_n)=
(x\mt z|\Delta(a)(\Om_n\mt y)).$$
And the third one:
$$(z|(\eta_{i(a^*a)}\mt id)W y)=(\Om_i\mt z|W(\hat{h}(a^*a)\mt y))=
(\Om_i\mt z|W(a^*\mt I)(\hat{h}(a)\mt y))=$$
$$=(\Om_i\mt z|\Delta(a^*)W(\hat{h}(a)\mt t))=
(\Delta(a)(\Om_i\mt z)|W(\hat{h}(a)\mt y))=$$
$$=(\hat{h}(a)\mt z|W^*(p_i\mt I)W\hat{h}(a)\mt y)
=(z|((\eta_a\mt id)W^*(p_i\mt I)W)y).$$

Now, let $\varphi$ be a positive linear functional on $A$ and 
$(K,\pi_\varphi,\Omega_\varphi)$ the associated GNS triple. We compute:
$$\eta_i((id\mt\varphi)\Delta(a^*a))=\varphi((\eta_i\mt id)\Delta(a^*a))=
\varphi((\eta_{i(a^*a)}\mt id)W)=\varphi((\eta_a\mt id)(W^*(p_i\mt I)W))$$

But since $W\in M(CB(L^2(G))\mt C^*_r(G_A))$
and $\sum_{i=1}^np_i$ converges strictly
(in $B(H)=M(CB(H))$) to $I$, we conclude that 
$\sum_i\eta_i((id\mt\varphi)\Delta(a^*a))=h(a^*a)\varphi(I)$. This 
shows  that $h((id\mt\varphi)\Delta(a^*a))$ is finite and 
equal to  $h(a^*a)\varphi(I)$. In this way we have proved:
\begin{prop}
Let $a\in D(\hat{h})$ and $\varphi$ be a positive linear functional on $A$. 
Then $$h((id\mt\varphi)\Delta(a^*a))=h(a^*a)\varphi(I).$$
\end{prop}
\dow\\
Therefore the weight $h$ is right-invariant.


\section{Putting all together.}

In this section we lift some objects from the section 5 to a
$C^*$-algebraic level and recapitulate the structure we got. 
Let us start by showing that $R$ defines $*$-antiautomorphism 
of $C^*_r(G_A)$. To this end it is enough to prove its continuity. 
This will immediately follow from
\begin{lem}
There exists an antiunitary operator $\hat{J}$ on $L^2(G_A)$ such that\vsm\vsm
$$\pi_{id}(R \om)=\hat{J}\pi_{id}(\om^*)\hat{J}\,\,,\,\om\in\sA(G_A).$$
\vsm\vsm
\end{lem}\vsm\vsm
{\em Proof: } Since $s_B$ (inverse of the groupoid $G_B$) is a diffeomorphism
of $G_A$, it defines a unitary operator on $L^2(G_A)$ (by a push-forward
of half-densities), which we also denote by $s_B$. On $L^2(G_A)$ there is
a canonical antiunitary involution, namely a complex conjugation: 
$\,\bar{}\,\,$.
Let us define  $\hat{J}:=\bar{}\cdot s_B$ and check that this is the 
right choice.
We choose $\psi_0=\ro\mt \nu$ --  a non-vanishing, real half-density on $G_A$  
as before lemma \ref{W'}. Then short computations show that:\notka{sB}
\begin{equation}\label{sB}
s_B(Yg\wedge gX)=\Adab(\al(g))^{-1}\Adaa(\ar(g))(-1)^{\dim\, A}
(Ys_B(g)\wedge s_B(g) X)
\end{equation}
$${\rm for \,\,}\,\, X\in\lma T_eA\,,\,Y\in\lma T_eB$$
From this formula we deduce that 
$(s_B\psi_0)(g)=|\Adab(\al(g))|^{1/2}|\Adaa(\ar(g))|^{1/2}\psi_0(g).$ 
Defining $(s_B\psi)$ for any smooth, compactly supported function $\psi$ on $G_A$
by $s_B(\psi\psi_0)=:(s_B\psi)\psi_0$ we get  
$$(s_B\psi)(g)=|\Adab(\al(g))|^{1/2}|\Adaa(\ar(g))|^{1/2}\psi(s_B(g)).$$ 
As usual, writing
$\om=f\om_0$,  we have to show that
$(Rf)*\psi=\overline{s_B((f^*)*\overline{(s_B\psi)})}$. We compute the left hand
side using the formulae for $R$ (\ref{defR}) and the multiplication 
(\ref{product-Ga}):
$$((Rf)*\psi)(g)=\int_B\mu_L^2(b)(Rf)(\al(g)b)\psi(\bl(\al(g)b)^{-1}g)=$$
$$=|\Adab(\al(g))|^{1/2}\int_B\mu_L^2(b)|\Adab(\ar(\al(g)b))|^{1/2} 
f(b^{-1}\al(g)^{-1})\psi(\bl(\al(g)b)^{-1}g).$$
And the right hand side
$$\overline{s_B(f^**\overline{s_B\psi})}(g)=
|\Adab(\al(g))|^{1/2}|\Adaa(\ar(g))|^{1/2} 
\overline{(f^**\overline{s_B\psi})(s_B(g))}=$$
$$=|\Adab(\al(g))|^{1/2}|\Adaa(\ar(g))|^{1/2}\int_B\mu_L^2(b)
\overline{(f^*)}(\al(s_B(g))b)(s_B\psi)(\bl(\al(s_B(g))b)^{-1}s_B(g))=$$
$$=|\Adab(\al(g))|^{1/2}|\Adaa(\ar(g))|^{1/2}\int_B\mu_L^2(b)
f(s_A(\al(g)^{-1}b))(s_B\psi)(\bl(\al(g)^{-1}b)^{-1}s_B(g))=$$
$$=|\Adab(\al(g))|^{1/2}|\Adaa(\ar(g))|^{1/2}\int_B\mu_L^2(b)
f(\bl(\al(g)^{-1}b)^{-1}\al(g)^{-1})(s_B\psi)(\bl(\al(g)^{-1}b)^{-1}s_B(g)).$$
Now apply  to the integral over $B$ a diffeomorphism 
$B\ni b\mapsto \bl(\al(g)^{-1}b)\in B$ to convert it into
$$|\Adab(\al(g))|^{1/2}|\Adaa(\ar(g))|^{1/2}\int_B\mu_L^2(b)
|\Adab(\ar(\al(g)b))| f(b^{-1}\al(g)^{-1})(s_B\psi)(b^{-1}s_B(g))$$
and use the definition of $s_B\psi$ to conclude that
that this is equal to left hand side.\\\dowl\\
Now we pass to the group $\tau_t$. The following proposition was proved in \cite{PS2}:
\begin{prop} Let $\G$ be a differential groupoid and $\sigma:\G\lra ]0,\infty[$
a smooth cocycle. The mapping 
$\sA(\G)\ni \om\mapsto \sigma_t(\om):= \sigma^{it}\om\in\sA(\G)$ extends to
a strongly continuous one parameter group on $C^*_r(\G)$. Let 
$\sigma_{\frac{i}{2}}$ be its analytic generator, then $\sA(\G)$ is a core for
$\sigma_{\frac{i}{2}}$ and  $\sigma_{\frac{i}{2}}(\om)=\sigma^{-1/2}\om$ 
for $\om\in\sA(\G)$.
\end{prop}
\dow

\noindent
From this proposition we infer that   $\tau_t$ defined by (\ref{deftau}) 
extends to a strongly continuous one parameter group on $C^*_r(G_A)$, moreover 
from prop. \ref{kap-tau-R} this group commutes with $R$. Since now we know
that $R$ $\tau_t$ and $\sigma_t$ are  continuous, using  lemma 
\ref{delta-R-tau} and prop. \ref{proteza} we easily obtain  the 
following equalities on $C^*_r(G_A)$:
$$\Delta\tau_t=(\tau_t\mt\tau_t)\Delta\,\,,\,\,\,\,
\Delta R=\sim(R\mt R)\Delta\,\,,\,\,\,\,
\Delta\sigma_t=(\sigma_t\mt\tau_t)\Delta.$$\vs

We finish this work with a short resume.
Let $(G;A,B)$ be a double Lie group. 
\begin{itemize}
\item There are naturally defined differential
groupoids $G_A, G_B$ over $A$ and $B$ respectively. 
Let $C:=C^*_r(G_A)$ be the reduced $C^*$-algebra of $G_A$.
\item The relation $\d:=m_B^T$ defines $\Delta\in Mor(C,C\mt C)$ which is
a comultiplication in the sense of the theory of quantum groups.
\item The inverse of the group $G$ defines the  involutive $*$-antiautomorphism
$R$ of $C$. 
\item There is a strongly continuous one parameter group of
$*$-automorphisms $\tau_t$ on $C$, which commutes with  $R$.
\item There exists right-invariant, densely defined, lower semi-contionuous
weight $h$ on $C$ which is a KMS weight with a modular group $\sigma_t$.
Moreover the weight $h$ is $\tau$ invariant.
\item The groups $\sigma_t$ and $\tau_t$ commute.
\item $\Delta\tau_t=(\tau_t\mt\tau_t)\Delta\,\,,\,\,\,\,
\Delta R=\sim(R\mt R)\Delta\,\,,\,\,\,\,
\Delta\sigma_t=(\sigma_t\mt\tau_t)\Delta.$
\end{itemize}

Because positions of groups $A$ and $B$ in  a DLG 
$(G;A,B)$ are completely symmetric, in fact, there are two quantum groups 
based on $C^*_r(G_A)$ and $C^*_r(G_B)$ respectively. 
One can think of them as duals, but the precise sense of this duality
is to be understood. If these algebras coincide with algebras defined by
the multiplicative operator $W$ in the ``standard way'' (see Appendix B), then
one can say that this is the meaning of duality. However, here we have more,
namely, there is a natural duality (in fact two ones) between $\sA(G_A)$ and
$\sA(G_B)$ (see Appendix A). This duality enable us to think about these 
algebras as ``dual Hopf algebras''. There is also a class of representations 
of (say) $C^*_r(G_B)$, the ones coming from morphisms  of differential groupoids 
$G_B \rel \G$, which define, via $W$, representations of quantum group 
based on $C^*_r(G_A)$.


\section{Appendixes}
{\large \bf A. Geometric interpretation of the function $Q$.}
Let $V$ be a finite dimensional vector space (over $C$ or $R$). Suppose
we are given four subspaces $L_1,L_2,R_1,R_2\subset V$, such that
$\dim L_1=\dim R_1$, $\dim L_2=\dim R_2$ and 
$V=L_1\ms R_2=R_1\ms L_2=L_1\ms L_2=R_1\ms R_2$.
Moreover let $\lambda_i,\rho_i\,,\,i=1,2$ be half-densities 
on $L_i$ and $R_i$  respectively. 
The quadruple $(\lambda_1,\rho_1,\lambda_2,\rho_2)$ defines 
two densities $d_1,d_2$ on $V$ as follows: 
$d_1:=(\lambda_1\mt\rho_2)(\lambda_2\mt\rho_1)$ and 
$d_2:=(\lambda_1\mt\lambda_2)(\rho_1\mt\rho_2)$. One can check, that 
$d_1,d_2$ depend only on $\lambda_1\mt\rho_1$ and $\lambda_2\mt\rho_2$ and the
dependence is bilinear, so in fact $d_1,d_2$ are bilinear mappings
$d_1,d_2:(\Om^{1/2}(L_1)\mt\Om^{1/2}(R_1))\times (\Om^{1/2}(L_2)\mt\Om^{1/2}(R_2))
\lra \Om^1(V)$. 
Since the space $\Om^1(V)$ is one dimensional we have that
$d_1=c d_2$ for some $c\in C$. To find constant $c$ is enough to compare
$d_1$ and $d_2$ on some basis in $V$. 
Let us recall that a $p$-density $d$ on $V$ is a mapping  
$d: V^n:=V\times\dots\times V\lra C\,,\,n=\dim V$ 
which satisfy the condition
$$d(v_{i_1}A_{i_11},\dots,v_{i_n}A_{i_nn})=|\det A|^p d(v_1,\dots,v_n)\,\,
{\rm (summation)},$$
for any $A\in M_{n\times n}(R)$ (or $A\in M_{n\times n}(C)$) and any 
$v_1,\dots,v_n\in V\,\,$\\
If $(v_1,\dots,v_n)=:{\bf v}$ is a basis in $V$ then we write the above
condition as $d({\bf v}A)=|\det A|^p d({\bf v})$.

\newcommand{\bfl}{{\bf l}}
\newcommand{\bfr}{{\bf r}}
\newcommand{\bftl}{{\bf \ti{l}}}
\newcommand{\bftr}{{\bf \ti{r}}}
\newcommand{\bfX}{{\bf X}}
\newcommand{\bfY}{{\bf Y}}
Choose ${\bf l}:=(l_1,\dots,l_n)$, $\,{\bf r}:=(r_1,\dots,r_n)$, $\,
{\bf \ti{l}}:=(\ti{l}_1,\dots,\ti{l}_m)$, $\,
{\bf \ti{r}}:=(\ti{r}_1,\dots,\ti{r}_m)$ -- bases in $L_1,R_1$, $L_2,R_2$, 
respectively. Let $I$ denotes the identity matrix 
(of appropriate dimension) and 
let matrices $A,B,C,D,J$,$K,G,H,$ $M,M_1,M_2$ be defined by:
$$(\bfl,\bftl)=:(\bfr,\bftr)\left(\begin{array}{cc} A & C\\B&D
\end{array}\right)
=:(\bfr,\bftr)M$$
$$(\bfl,\bftl)=:(\bfr,\bftl)\left(\begin{array}{cc} J & 0\\K&I
\end{array}\right)
=(\bfr,\bftl)M_1
\,\,,\,\,\,
(\bfl,\bftl)=:(\bfl,\bftr)\left(\begin{array}{cc} I& G\\0&H
\end{array}\right)=(\bfl,\bftr)M_2.$$
Now we can compare $d_1$ and $d_2$
$$d_1(\bfl,\bftl):=(\lambda_1\mt\rho_2)(\bfl,\bftl)(\lambda_2\mt\rho_1)(\bfl,\bftl)
=(\lambda_1\mt\rho_2)((\bfl,\bftr)M_2)(\lambda_2\mt\rho_1)((\bfr,\bftl)M_1)=$$
$$=|\det H|^{1/2}|\det J|^{1/2}\lambda_1(\bfl)\rho_2(\bftr)
\lambda_2(\bftl)\rho_1(\bfr).$$
$$d_2(\bfl,\bftl):=(\lambda_1\mt\lambda_2)(\bfl,\bftl)(\rho_1\mt\rho_2)(\bfl,\bftl)=
\lambda_1(\bfl)\lambda_2(\bftl)(\rho_1\mt\rho_2)((\bfr,\bftr)M)=$$
$$=|\det M|^{1/2}\lambda_1(\bfl)\lambda_2(\bftl)\rho_1(\bfr)\rho_2(\bftr).$$
In this way 
$\,d_2=\frac{|\det M|^{1/2}}{|\det H|^{1/2}|\det J|^{1/2}} d_1$.\\
But from equalities: 
$$\bfl=\bfr A+\bftr B\,,\,
\bftl=\bfr C +\bftr D\,,\, \bfl=\bfr J+\bftl K\,,\,\bftl=\bfl G+\bftr H,$$
we infer that 
$C=AG\,,\,B=DK\,,\,J=A(I-GK)\,,\,D=(I-KG)$ and we can write 
$$d_2=\left|\frac{\det A\det D}{\det M} \right|^{1/2} d_1.$$

Now for a DLG  $(G;A,B)$ let us fix a point $g\in G$ and define
$$L_1:=T^{l}_gG_B\,,\,L_2=T^{l}_gG_A\,,\,R_1:=T^{r}_gG_B\,,\,R_2:=T^{r}_gG_A
{\rm\,\, and\,\,} V:=T_gG.$$
Let $\bfX:=(X_1,\dots,X_n)\,,\,\bfY=(Y_1,\dots,Y_m)$ be bases 
in $T_eA$ and $T_eB$, respectively. Then $g\bfX,g\bfY$, $\bfX g,\bfY g$ 
are bases in $L_1,L_2,R_1,R_2$, respectively. Moreover, if 
$M=\left(\begin{array}{cc} A & C\\B& D\end{array}\right)$ is a matrix of 
adjoint representation of $g$ in $(\bfX,\bfY)$, then 
$(g\bfX,g\bfY)=(\bfX g, \bfY g)M$.
Comparing this with the definition of the function $Q$ we conclude that
$d_1=|Q(g)|^{1/2} d_2$.

In this way we see that the function $Q$ relates two natural
dualities between $\sA(G_A)$ and $\sA(G_B)$ given, for 
$\om_A=\lambda_A\mt \rho_A\in\sA(G_A)$ and  $\om_B=\lambda_B\mt\rho_B\in \sA(G_B)$, 
by formulae:
$$<\om_A,\om_B>_1:=\int_G (\lambda_A\mt\rho_B)(\lambda_B\mt\rho_A)\,\,
{\rm and}\,\,
<\om_A,\om_B>_2:=\int_G (\lambda_A\mt\lambda_B)(\rho_A\mt\rho_B).$$
{\large\bf B. Comparison with the standard approach.}
Let $W\in B(H\mt H)$ be a multiplicative unitary operator.
As shown by Baaj and Skandalis \cite{BajSka} the set:
$$C:=\overline{\{(\eta\mt id)W\,,\,\eta\in B(H)_*\}}=
\overline{span\{(\eta_{xy}\mt id)W\,,\,x,y\in H\}},$$
(bar denotes norm closure) is a $C^*$-algebra with a comultiplication.
This is the ``standard procedure''.
In this appendix we consider the following problem. 
For a DLG $(G;A,B)$  there is a multiplicative unitary (manageable) 
$W\in B(L^2(G)\mt L^2(G))$ and we have 
$C^*$-algebra $C$ defined as above. It is easy to see that
in the definition of $C$, it is enough to consider vectors $x,y$
from a dense subspace of $H$.
So a ``typical'' element of $C$ is obtained from two vectors 
$x,y$ by $a_{xy}:=(\eta_{xy}\mt id)W$, in  particular it is defined by 
two smooth half-densities on $G$  with compact support.
On the other hand we have the groupoid $G_A$ and its reduced $C^*$-algebra
defined in a different way, and its ``typical'' elements are bidensities 
on $G_A$. We can suppose that those algebras are in fact 
the same. We will not prove this equality here, but we will show that
$C\subset C^*_r(G_A)$. From the remark we made above, it follows that
we have to interpret elements $(\eta_{xy}\mt id)W$ as elements of $\sA(G_A)$ for 
$x,y$ --  smooth with compact support. The question is: how from two smooth 
half-densities on $G$ and $W$ we can get an element of $\sA(G_A)$ ? 

Suppose we are given $\varphi,\psi$ -- half-densities on $G$ with compact support. 
Let $U$ be the bisection implementing $W$ and let $g\in G$ be given. 
Consider the set $U (g, A)\in G\times G$.
It is easy to see that $U (g, A)=\d(g)$ (we recall: $\d=m_B^T$). 
Let $(g_1,g_2)\in \d(g)$. There are the  following natural isomorphisms:
\begin{eqnarray}
 \Om^{1/2}(T_{g_1}G) & {\simeq} & 
\Om^{1/2}(T_{\bl(g)}B)\mt\Om^{1/2}(T^l_{g_1}G_B)  \nonumber\\
  \Om^{1/2}(T_{g_2}G) & \simeq & 
\Om^{1/2}(T_{\br(g)}B)\mt\Om^{1/2}(T^r_{g_2}G_B) \nonumber\\
 \nonumber \Om^{1/2}(T_{\bl(g)}B) & \simeq & \Om^{1/2}(T^r_gG_A)\\
 \nonumber \Om^{1/2}(T_{\br(g)}B) & \simeq  & \Om^{1/2}(T^l_gG_A)\\
\nonumber \Om^{1/2}(T_{(g_1,g_2)}\d(g)) & \simeq & \Om^{1/2}(T^l_{g_1}G_B)\\
\nonumber \Om^{1/2}(T_{(g_1,g_2)}\d(g)) & \simeq & \Om^{1/2}(T^r_{g_2}G_B)
\end{eqnarray}
The first two isomorphisms follow  from the fact that $\bl$ and $\br$ are 
surjective submersions;\\
The third and the forth ones are given, respectively, by ismorphisms 
$T_{\bl(g)}B\ni X\mapsto X\ar(g) \in T^r_gG_A$ and 
$T_{\br(g)}B\ni Y\mapsto \al(g) Y \in T^l_gG_A$;\\
The last two ones are due to the fact that 
$\pi_1: \d(g)\ni(g_1,g_2)\mapsto g_1 \in F^l_B(g)$
and $\pi_2: \d(g)\ni(g_1,g_2)\mapsto 
g_2 \in F^r_B(g)$ are  diffeomorphisms.\\
And consequently:
\begin{eqnarray}
\Om^{1/2}(T_{g_1}G) & \simeq & \Om^{1/2}(T^r_gG_A) \mt 
\Om^{1/2}(T_{(g_1,g_2)}\d(g))\nonumber\\
\Om^{1/2}(T_{g_2}G) & \simeq & \Om^{1/2}(T^l_gG_A)\mt 
\Om^{1/2}(T_{(g_1,g_2)}\d(g)).\nonumber
\end{eqnarray}

In this way, having $\varphi(g_1),\psi(g_2)$ we can define a 1-density on 
$T_{(g_1,g_2)}\d(g)$ with values in one dimensional vector space 
$\Om^{1/2}(T^l_gG_A)\mt\Om^{1/2}(T^r_gG_A)$. Let us denote the resulting 
mapping by $\Phi_g$.
The explicit formula for $\Phi_g(\varphi,\psi)$ is given by
$$[\Phi_g(\varphi,\psi)(u)](v\wedge w)=\varphi(g_1)(wa\wedge \pi_1(u))
\psi(g_2)(\ti{a} v\wedge \pi_2(u)),{\rm where}$$ 
$$u\in\lma T_{(g_1,g_2)}\d(g)\,,\,v\in\lma T_g^lG_A\,,\,w\in\lma T_g^rG_A
\,,\,\,a:=g^{-1}g_1\,,\,\,\ti{a}:=g_2g^{-1}$$
 and  $\pi_1\,,\,\pi_2$ are defined above.\\
Integrating $\Phi_g(\varphi,\psi)$ we get
a bidensity on $G_A$, which will be denoted by $\Phi(\varphi,\psi)$. 
Finally,  we define a mapping
$$\hat{\Phi}: (\varphi,\psi)\mapsto 
|Q|^{1/2}\Phi(\overline{s_B(\varphi)},\psi) {\rm \,\,\,i.e.\,\,\,\,}
\hat{\Phi}(\varphi,\psi)(g)=
|Q|^{1/2}(g)\int_{\d(g)}\Phi_g(\overline{s_B(\varphi)},\psi).$$
and will prove  the following
\begin{prop} Let $\varphi,\psi$ be smooth, compactly supported 
half-densities on $G$. Then
$\hat{\Phi}(\varphi,\psi)\in\sA(G_A)$, 
$\kappa(\hat{\Phi}(\varphi,\psi))=(\hat{\Phi}(\psi,\varphi))^*$ 
and $(\eta_{\varphi\psi}\mt id)W=\pi_{id}(\hat{\Phi}(\varphi,\psi))$.
\end{prop}
{\em Proof: } We begin by computing the mapping $\hat{\Phi}$. 
Let us choose $\mu_0,\nu_0$ -- real, non-zero half-densities on 
$T_eB$ and $T_eA$, respectively, and let $\mu_l,\nu_r$ denote the corresponding
left and right-invariant half-densities on $B$ and $A$. 
Then $\nu_r\mt\mu_l$ is a real, non-vanishing half-density on $A\times B$ 
and, since the mapping $A\times B\ni(a,b)\mapsto ab \in G$ is a 
diffeomorphism, this half-density defines a half-density on $G$ which will be 
denoted by $\epsilon$. Explicitly
$$\epsilon(g)(Xg\wedge gY)=\nu_0(X)\mu_0(Y)\,,\,{\rm for } X\in \lmA\,,\,
Y\in \lmB.$$
Now we can write 
$\varphi=:f_\varphi\eps\,,\,\psi=:f_\psi\eps$ and 
$\hat{\Phi}(\varphi,\psi)=:f_{\varphi\psi}\om_0$.
Let us define a mapping
$$\alpha: A\ni a\mapsto
\alpha(a):=(m_B(g,s_B(ag)),ag)=(g\ar(ag)^{-1},ag)\in\d(g).$$
This is clearly a diffeomorphism.
Let us put $v:=g Y\,,\,w:=\ti{Y}g\,,\,Y,\ti{Y}\in\lma T_eB\,,\,
\ti{u}:=X a \,,\,X\in\lma T_eA$,
then $\pi_1\alpha(\ti{u})=
(-1)^{\dim\,A}\Adba(\bl(a\bl(g))^{-1} s_B(a\bl(g)) X$ 
and $\pi_2\alpha(\ti{u})=X a g$. 
Inserting this into the formula for $\Phi_g$ and
using (\ref{sB}) we get
$$(\Phi_g(\overline{s_B(\varphi)},\psi)(\alpha(Xa))(gY\wedge \ti{Y}g)=$$
$$=|\Adba(\bl(a\bl(g)))|^{-1/2}
\overline{s_B(\varphi)}(\ti{Y}s_B(a\bl(g))\wedge s_B(a\bl(g))X)\,
\psi(ag Y\wedge Xag)=$$
$$=\left|\frac{\Adab(a)}{\Adaa(\ar(a\bl(g)))\Adba(\bl(a\bl(g))) }\right|^{1/2}
\overline{\varphi(\ti{Y}a\bl(g)\wedge a\bl(g) X )}\,\psi(ag Y\wedge Xag),$$
and since
$$gX\wedge Yg=\frac{\Adaa(\ar(g))\Adba(\bl(g))}{\Adab(\al(g))\Adbb(\br(g))}
(X g\wedge g Y),$$
we obtain
$$\Phi_g(\overline{s_B(\varphi)},\psi)(\alpha(Xa))(gY\wedge \ti{Y}g)=
|\Adbb(\bl(g))|^{-1/2}
\overline{\varphi(a\bl(g)\ti{Y}\wedge X a\bl(g))}\,\psi(ag Y\wedge Xag)$$
Finally (we also use (\ref{Q-mod}))
\begin{equation}\label{f-fi-psi}
f_{\varphi\psi}(g)=\left|\frac{\Adab(\ar(g))}{\Adbb(\br(g))}\right|^{1/2}
\int_A \nu_r^2(a)\overline{f_\varphi(a\bl(g))}f_\psi(ag).
\end{equation}
From this equation is clear that $f_{\varphi\psi}\om_0\in\sA(G_A)$.\vs

Having the above expression  and using (\ref{def-kappa}) one easily
proves that $\kappa (\hat{\Phi}(\varphi,\psi))=\hat{\Phi}(\psi,\varphi)^*$.\vs

Now we are going to prove the last equality. 
Of course it is enough to prove that 
$(\varphi\mt z|W(\psi\mt t))=(z|\pi_{id}(\hat{\Phi}(\varphi,\psi))t)$ for 
compactly supported, smooth half-densities $z,t$ on $G$. 
Let  $\psi_0=:\ro\mt\nu_l$, where $\nu_l$ is the left-invariant half-density on
$A$ defined by $\nu_0$. (this is $\psi_0$ as defined before the proof
of lemma \ref{W'}). We can write  $t=:f_t\psi_0$ and $z=:f_z\psi_0$.

Using the formula (\ref{product-Ga}) we have 
$\pi_{id}(\hat{\Phi}(\varphi,\psi))t=(f_{\varphi\psi}*f_t)\psi_0$ and
$$(f_{\varphi\psi}*f_t)(g)=\int_B \mu_l^2(b)f_{\varphi\psi}(\al(g)b)
f_t(\bl(\al(g)b)^{-1}g)=$$
$$=\int_B\mu_l^2(b)|\Adab(\ar(\al(g)^{-1}b))|
f_{\varphi\psi}(\al(g)\bl(\al(g)^{-1}b)
f_t(b^{-1}g)=$$
$$=\int_B\mu_l^2(b)|\Adab(\ar(\al(g)^{-1}b))|f_t(b^{-1}g)
\left|\frac{\Adab(\ar(\al(g)\bl(\al(g)^{-1}b)))}
{\Adbb(\bl(\al(g)^{-1}b))}\right|^{1/2}\times$$
$$\times 
\int_A\nu_r^2(a)\overline{f_\varphi(ab)} f_\psi(a \al(g)\bl(\al(g)^{-1}b))=$$
$$\int_G \epsilon^2(\ti{g})\left|\frac{\Adab(\ar(\al(g)^{-1}\br(\ti{g})))}
{\Adbb(\bl(\al(g)^{-1}\br(\ti{g}))}\right|^{1/2}
\overline{f_\varphi(\ti{g})}f_\psi(\ti{g}\al(\br(\ti{g})^{-1} g))
f_t(\br(\ti{g})^{-1}g)=$$
$$\int_G \epsilon^2(\ti{g})\left|\frac{\Adab(\ar(\al(g)^{-1}\br(\ti{g})))}
{\Adbb(\bl(\al(g)^{-1}\br(\ti{g}))}\right|^{1/2}
\overline{f_\varphi(\ti{g})}\,(f_\psi\mt f_t)(W^{-1}(\ti{g},g))$$
In the second equality we use a diffeomorphism $B\ni b\mapsto \bl(\al(g)b)\in B$,
then the expression for $f_{\varphi\psi}$ and finally the definition of $\epsilon$ 
($\ti{g}:=ab$), and formula for $W^{-1}$ (\ref{W}).

On the other hand, doing computations similar to the ones done in the 
proof of lemma \ref{W-h}  we get
$(\varphi\mt z|W(\psi\mt t))=(z|y)$ for $y=:f_y\psi_0$ and $f_y$ given by
$$f_y(g)=\int_G \epsilon^2(\ti{g})
\overline{f_\varphi(\ti{g})}[W(f_\psi\mt f_t)](\ti{g}, g),$$
where  $W(f_\psi\mt f_t )$ is defined by
$W(f_\psi\epsilon\mt f_t\psi_0)=:W( f_\psi\mt f_t)(\epsilon\mt\psi_0).$

Since we can easily compare $\epsilon\mt\psi_0$ with $\psi_0\mt\psi_0$,
the expression for the function $W(f_t\mt f_\psi)$ follows from
the action of $W$ on $\psi_0\mt\psi_0$ given in the proof of
lemma \ref{W'}. In this way we obtain
$$W(f_\psi\mt f_t )(\ti{g},g)=(f_\psi\mt f_t )(W^{-1}(\ti{g},g))
\left|\frac{\Adab(\ar(\al(g)^{-1}\br(\ti{g})))}
{\Adbb(\bl(\al(g)^{-1}\br(\ti{g}))}\right|^{1/2}$$
Now insert this  into formula for $f_y$ and compare with
the expression for $f_{\varphi\psi}*f_t$. This completes the proof
of the proposition.\\\dowl


\ed